\documentclass[english]{article}
\usepackage[T1]{fontenc}
\usepackage{amsmath}
\usepackage{amssymb}
\usepackage{amsthm}
\usepackage{array}
\usepackage{epsfig}
\usepackage{graphicx}
\usepackage{babel}
\usepackage{algorithm}
\usepackage[noend]{algpseudocode}

\usepackage{authblk}

\usepackage{etoolbox}
\usepackage{tikz}
\usetikzlibrary{tikzmark}
\usetikzlibrary{calc}

\errorcontextlines\maxdimen

\newcommand{\ALGtikzmarkcolor}{black}
\newcommand{\ALGtikzmarkextraindent}{4pt}
\newcommand{\ALGtikzmarkverticaloffsetstart}{-.5ex}
\newcommand{\ALGtikzmarkverticaloffsetend}{-.5ex}
\makeatletter
\newcounter{ALG@tikzmark@tempcnta}

\newcommand\ALG@tikzmark@start{%
    \global\let\ALG@tikzmark@last\ALG@tikzmark@starttext%
    \expandafter\edef\csname ALG@tikzmark@\theALG@nested\endcsname{\theALG@tikzmark@tempcnta}%
    \tikzmark{ALG@tikzmark@start@\csname ALG@tikzmark@\theALG@nested\endcsname}%
    \addtocounter{ALG@tikzmark@tempcnta}{1}%
}

\def\ALG@tikzmark@starttext{start}
\newcommand\ALG@tikzmark@end{%
    \ifx\ALG@tikzmark@last\ALG@tikzmark@starttext
    \else
        \tikzmark{ALG@tikzmark@end@\csname ALG@tikzmark@\theALG@nested\endcsname}%
        \tikz[overlay,remember picture] \draw[\ALGtikzmarkcolor] let \p{S}=($(pic cs:ALG@tikzmark@start@\csname ALG@tikzmark@\theALG@nested\endcsname)+(\ALGtikzmarkextraindent,\ALGtikzmarkverticaloffsetstart)$), \p{E}=($(pic cs:ALG@tikzmark@end@\csname ALG@tikzmark@\theALG@nested\endcsname)+(\ALGtikzmarkextraindent,\ALGtikzmarkverticaloffsetend)$) in (\x{S},\y{S})--(\x{S},\y{E});%
    \fi
    \gdef\ALG@tikzmark@last{end}%
}

\apptocmd{\ALG@beginblock}{\ALG@tikzmark@start}{}{\errmessage{failed to patch}}
\pretocmd{\ALG@endblock}{\ALG@tikzmark@end}{}{\errmessage{failed to patch}}
\makeatother

\usepackage{xy}

\usepackage[draft, textwidth=0.85in]{todonotes}

\newcommand{\ccdot}{\;\cdot\;}

\usepackage{setspace}
\usepackage{esint}

\newcommand{\Idx}{\mathop{\mathrm{Idx}}}
\newcommand{\Law}{\mathop{\mathrm{Law}}}

\newtheorem{theorem}{Theorem}[section]

\newtheorem{lemma}[theorem]{Lemma}
\newtheorem{remark}[theorem]{Remark}
\newtheorem*{assumption}{Assumption}

\newlength{\bibitemsep}
\newlength{\bibparskip}
\let\oldthebibliography\thebibliography
\renewcommand\thebibliography[1]{%
  \oldthebibliography{#1}%
  \setlength{\parskip}{\bibitemsep}%
  \setlength{\itemsep}{\bibparskip}%
}

\author[$\star$]{Gabriel Earle}
\author[$\dagger$]{Jonathan Mattingly}
\affil[$\star$]{Department of Mathematics, University of Massachusetts-Amherst}
\affil[$\dagger$]{Departments of Mathematics and of Statistical 
  Science, Duke University}
\title{{Convergence of Stratified MCMC Sampling of Non-Reversible Dynamics}}
\usepackage{geometry}
\algnewcommand\algorithmicto{\textbf{to}}
\begin{document}

\maketitle

\abstract We present a form of stratified MCMC algorithm built with
non-reversible stochastic dynamics in mind. It can also be viewed
as a generalization of the exact milestoning method, or form of
NEUS.
We prove convergence of the method under certain assumptions, with
expressions for the convergence rate in terms of the process's behavior
within each stratum and large scale behavior between strata. We show
that the algorithm has a unique fixed point which corresponds to the
invariant measure of the process without stratification. We will show
how the speeds of two versions of the new algorithm, one with an extra
eigenvalue problem step and one without, relate to the mixing rate
of a discrete process on the strata, and the mixing probability of
the process being sampled within each stratum. The eigenvalue problem
version also relates to local and global perturbation results of discrete
Markov chains, such as those given by Van Koten, Weare et. al.

\section{Introduction} \label{introduction}

Markov Chain
Monte Carlo (MCMC) is an often used method to produce samples from a distribution, when a Markov kernel converging 
to that distribution is known.
Stratification of  MCMC methods is a well-studied
form of rare event sampling. Cases of interest include systems where
regimes of low probability have outsize importance, or systems with
multiple regimes of high probability but rare transitions between
them. In such cases, the sample space can be broken up into smaller
parts, called strata, and a Markov Chain can be run within each. The
results are then combined in some way, allowing for an estimate of
a distribution on the whole space which can be obtained more quickly
than via a simple MCMC method. There are many advantages of dividing
the pace into strata. Including the fact that many of the computations
can be run in parallel and computational resolution and effort can be
concentrated in regions of interest. 

In this article, we present and prove convergence results for a specific
stratified MCMC scheme close to those in \cite{dinner_mattingly_tempkin_koten_weare_2018}. The algorithm, which we call the ``injection
measure method'', constructs an estimate of the invariant measure
of a Markov chain, and is built specifically with non-reversible
Markov chains
in mind. It can also be viewed as a version of Non-Equilibrium Umbrella
Sampling (NEUS), as in \cite{dinner_mattingly_tempkin_koten_weare_2018},
or as an extension of the exact milestoning method, as detailed in
\cite{aristoff_bello-rivas_elber_2016}. Non-reversible Markov chains
sampling problems typically arise in two settings. First, when one
samples an invariant measure which does not satisfy
detailed balance; and  hence, has a nontrivial flux through the system
in its stationary state. Such stationary states are often referred to as
non-equilibrium steady-states. Second, when one samples the  non-reversible Markov chain
can be obtained by adding time to the dynamics as one of the
state-variables. Because of direction of time ensures that the system
has a non trivial flux. As  illustrated in
\cite{dinner_mattingly_tempkin_koten_weare_2018}, the resulting space-time dynamics can be used to
study out-of-equilibrium  transitions rates and other transient
phenomenon. See \cite{dinner_mattingly_tempkin_koten_weare_2018} for
more details.

The key object of the method, and our analysis of it, is a collection
of distributions on the strata known as injection measures. These
estimate how particles following the Markov chain are likely to be
distributed on the step when they enter one stratum from another.
Also associated to each stratum is a corresponding weight, estimating
how likely particles are to enter the stratum, relative to the other
strata. If the injection measures and weights were known exactly,
they could be used to compute estimates of the invariant measure within
each stratum, and patch those estimates together with the correct
weights. In practice, the injection measures and weights will most
likely not be known, and so they must themselves be estimated. We
propose estimating them iteratively, finding new injection measures
via trajectories started from the current ones. This iteration is
the main step of our formulation and those discussed in \cite{dinner_mattingly_tempkin_koten_weare_2018}.

The corresponding weights can be calculated in two ways, leading to
two versions of the injection measure method. In the first, the weights
are found by applying a transition matrix determined by the measures
to the previous weights. We call this form of the method the basic
version. In the other, the principle eigenvector of the matrix is
found, and its entries are taken to be the new weights. This is called
the eigenvector version. The eigenvector version performs strictly
better in our numerical experiments, but we as yet only have a proof
of local convergence for it, whereas we prove global convergence of
the basic version.

In order to prove our main results, we need some assumptions about
the behavior of the Markov chain being sampled and how it interacts
with the chosen strata. In addition to some standard regularity assumptions,
we will need to assumptions on two types of behavior. The first is
a ``microscopic'' assumption, governing how particles following
the chain move within a stratum. In effect, our assumption will be
that any two particles starting somewhere in the same stratum have
some chance of being coupled at or before the time of exit. We will
also need a ``macroscopic'' assumption, about how the chain moves
overall mass between strata. We will have several forms of this, depending
on the version of the algorithm and precision of the theorem we wish
to prove. However, each form of the assumption roughly states that
the chain moves mass between the various strata at a suitable rate,
or with enough regularity. The algorithm's overall convergence speed
can then be expressed in terms of the ``microscopic'' coupling probability
and the ``macroscopic'' rate. That is the substance of our main
theorems.

Specifically, our first result shows that, under the assumptions,
the algorithm has a unique fixed point, and the injection measures
given by the fixed point are correspond to the original invariant
measure. Next, we prove that, if $c$ is the microscopic coupling
probability and $\lambda$ is the macroscopic rate of convergence,
then the algorithm converges to the fixed point in total variation
at a rate that roughly looks like 
\[
O(\max(\sqrt{1-c},\lambda)).
\]
We also show that, for a constant $r$ reflecting the sensitivity
of the weights to the entries of a macroscopic transition matrix,
the eigenvector version converges, for sufficiently good starting
estimates, at a rate approaching 
\[
1-\frac{c}{1+r}.
\]
More precisely, we will show that, if $G$ is the matrix of transition rates between strata for the dynamics in equilibrium, 
then $r$ bounds the relative sensitivities of the invariant distribution of $G$ to small changes in the entries of $G$. 

The structure of this article will be as follows. In Section~\ref{inj_meas_method}, we
give the basic notation needed for our results, and state the two
versions of the injection measure method explicitly. Section~\ref{results},
we outline the assumptions needed and state the main theorems we can
now prove under them. In Section~\ref{numerics}, we give some results of numerical
simulations of the method on a simple example system, and relate them
to some of our theoretical results. In Section~\ref{proofs}, we give the proofs
of our theorems in detail. 

\section{The Injection Measure Method} \label{inj_meas_method}

\subsection{Intuition Behind the Algorithm}

To motivate the steps of the method, consider a space broken up into
$J$ subsets, or strata, and a process $X_{j}^{n}$,
following a kernel $P$ and confined to the $j$-th stratum. If we wish to sample the
invariant distribution restricted to that stratum, then the question becomes
what to do if a step drawn according to $P$ would have $X_{j}^{n}$
leave the stratum. One option would be to simply ``bounce it back''
from the boundary, which would mean setting $X_{j}^{n+1}=X_{j}^{n}$.
We could then keep track of how many attempted exits occur between
the strata and use that to form a matrix $G$ on the set of strata, or $\{1,\ldots,J\}$. The empirical distributions
of the $X_{j}^{n}$ can then be combined, with weights given by the
invariant distribution of $G$.

The above scheme works under a key assumption: that the kernel $P$
is reversible. The reason this is needed is that ``bouncing a particle
back'' when it tries to leave a stratum is like saying that another
particle comes to replace it at the exact location it left
from. The assumption that this occurs for many particles in equilibrium
is exactly the assumption that detailed balance holds, at least along
the strata boundaries. But we are precisely interested in the non-reversible
case. Therefore, we must have a better idea of how, and if, a particle
leaving  a stratum is replaced.

Thus we propose building an injection measure for each stratum, which captures how 
particles entering that stratum are distributed. In addition, we need weights that capture how many particles
enter each stratum, relative to the other strata. Of course, we cannot expect these to reflect how particles entering a stratum are in equilibrium, at first. Therefore, we will use the starting injection measures to calculate the distributions of particles leaving each stratum, called exit measures. Then, for each $j$, we look at the distribution of all particles leaving other strata, via their exit measures, and entering stratum $j$. We use this to form a new injection measure and weight for the $j$-th stratum, for all $j$. Thus we have an iterative procedure, building new injection measures and weights from old ones. We wish to show that the injection measures thus defined converge over time to some fixed point, called the equilibrium injection measures, and that they give us an estimate of the invariant measure of $P$. We call the method thus outlined the basic version of our algorithm. 

There is one further step we can add to the method. Just as in the reversible case, we can make a transition matrix $G$ based on how many particles leave each stratum into each other one, or how much mass the $k$-th exit measure gives ot the $j$-th stratum. Then we can replace the old weights on each iteration by the invariant vector of $G$. We call the method, with this added step, the eigenvector version of the injection measure method. We will be able to show, under certain assumptions, an improved bound for the rate of convergence for this version, with the drawback that we have not yet shown that it converges for any starting injection measures. That is, our theorems only show that it converges to the correct distribution locally. 

\subsection{Setting and Notation}

Suppose we are interested in sampling from the invariant distribution
$\pi$ of a Markov kernel $P,$ defined on a state space
$A\subset\mathbb{R}^{d}$.
The simplest version of the injection measure method of stratified
sampling is to break $A$ into subsets,
or strata, $A_{1},\ldots A_{J}$ which partion the space.
One then runs the
Markov chain $P$, starting
from a measure  $\nu_{j}^0$ concentrated  with in each strata $A_j$, until it exist. The algorithm calculates new starting
measures $\{\nu_{j}^1\}_j$ biased on the exists from the collection of
$\{A_j\}_j$. This process, which is   iterated until the  measure
$\{\nu_{j}^{n}\}_j$ converge
also provides a collection of weights $\{a_j\}$ and occupation measure
$\pi_j$ in each $A_j$ so that the desired
target measure $\pi= \sum a_j \pi_j$.

It the above setting is clear which of the strata $A_j$ contain the
process  at any given times since they are disjoint.
However, we are also
interested in the case where there is overlap between the strata,
and the point at which a particle is declared to have left one stratum
and entered another is possibly random. Though we will consider a more
general formulation, for the moment consider the following
illustrative setting.

Consider the following setup: we have strata $A_{1},\ldots, A_{J}$
covering $A$ as before, but the strata are no longer
disjoint. We also assume we have a partition of unity
$\psi_{1},\ldots,\psi_{J}$ whose suport coincides with the
$A_{1},\ldots, A_{J}$. Namely $\mathrm{supp}(\psi_j)=A_j$ and for all
$x \in A$, $\psi_j(x)
\geq 0$ and
\begin{align*}
  \sum_j \psi_j(x) = 1\,.
\end{align*}

When a particle enters
$A_{j}$at $X_{0}\sim\nu_{j}$, a value $\kappa$ is chosen from some
distribution $\eta$ on $[0,1]$. The particle then moves according
to $P$, with its position at time $n$ being $X_{n}\sim
P(X_{n-1},\ccdot)
$. When $\psi_{j}(X_{n})<\kappa$,
the particle is declared to have exited $A_{j}$. The index of the
stratum it exits into is chosen from $1,\ldots j-1,j+1,\ldots J$,
with probabilities proportional to by $\psi_{1}(X_{n}),\ldots,\psi_{j-1}(X_{n}),\psi_{j+1}(X_{n}),\ldots,\psi_{J}(X_{n}).$
We want our results to cover this more general setting, because it
is often computationally useful. In particular, one does not need to
insure that the $A_{1},\ldots, A_{J}$ are disjoint which can
be computationally intensive. Additionally, the ``softening'' of the
exit boundary by introducing random exit  times seems to soften any
artifacts from the strata edges.

We can fit this more  general case,  into our original setting by
consider the following augmented space 
\[
A^{\prime}=\{(x,k):x\in A_{k},\,1\leq k\leq J\}
\]
 with strata given by 
\[
A_{j}^{\prime}=\{(x,j):x\in A_{j}\}.
\]
The $A_{j}^{\prime}$ are then a partition of $A^{\prime}$ even when the $A_j$
are only a cover of $A$ (namely, not disjoint).

We then extend the Markov dynamics of $P$ to $A'$ in the stratified
setting by defining a collection of kernels $\{P'_\kappa :
\kappa \in [0,1]\}$ on $A'$ as follows. If the initial state is
$(X_0,j_0)$ then the new state $(X_1,j_1)$ is constructed as
follows: $X_1\sim P(X_0,\ccdot)$ and $j_1=j_0$ if $\psi_{j_0}(X_1)
\geq\kappa$. If $\psi_{j_0}(X_1)
< \kappa$ then $j_1$ is chosen randomly according the probabilities
$P(j_1 = k) \propto \psi_k(X_1)$ for $k \neq j_0$ as described above.  
Then, when a particle is started in $A_{j}^{\prime}$ via the injection
measure, a value of $\kappa$ is drawn according to some probability
measure $\eta$ on $[0,1]$ , and the particle moves by
$P_{\kappa}^{\prime}$ until it leaves $A_{j}^{\prime}$.

In light of this construction, we can recast this more setting into
the initial frame work of disjoint intervals.

Since all of the kernels
$P_{\kappa}^{\prime}$ have the same action on x, i.e. the x-marginal of
 $P_{\kappa}^{\prime}((x,j),\cdot)$ is independent of $\kappa$, and
the $x$-marginal of the invariant distribution of each kernel is
$\pi$. In particular, the kernel $P ^{\prime} (x,\cdot)=\int_{[0,1]}P_{\kappa}^{\prime}(x,\cdot)\eta(d\kappa)$ 
acts the same on $x$ as $P$. 

For the remainder of this article, we will use this last setting.
We will drop the superscript of $A^{\prime}$, $P_{\kappa}^{\prime}$, etc.
and simply assume that the strata $A_{j}$ are a partition of $A$,
and that we have a family of Markov kernels $P_{\kappa}$ on $A$.
We will also require that $P(x,\cdot)=\int_{[0,1]}P_{\kappa}(x,\cdot)\eta(d\kappa)$,
so integrating the kernels over $\kappa$ gives us back the kernel
of the original process we wanted to sample. We will show that, in
this setting, the injection measure method gives a way of approximating
the invariant distribution of $P$. In fact, our results apply to any family of kernels
$P_{\kappa}$ to an collection of kernels $P_{\kappa}$ which we will
make more explicit later. However, our primary interest  is the
specific 
choice of kernels outlined above. 

\begin{remark}
  Note that, if $\kappa$ is always chosen to be $1$, so that 
there is only one kernel $P_{1}$, then we are back in the original 
setting, so our results here cover all the cases in which we are interested. 
\end{remark}

Since the $A_{j}$ are now assumed to be disjoint, we can define the
index of a point in $A$ as follows: $\Idx(x)=j$ if $x\in A_{j}$. 

Now, we are interested in the point at which the sample process exits
a stratum, so, for any $X_0=x\in A,$ define the exit time starting from to be 
\begin{align*}
  \tau=\inf\Big\{ n\geq1:\Idx(X_{n})\neq\Idx(X_0)\Big\}
\end{align*}
where $\kappa$ is first chosen according to $\eta$ and then $\{X_{n}\}_{n\geq0}$
moves according to $P_{\kappa}$. Note that, since $\Idx(X_{0})$
is possibly random, we define the exit time starting from $n=1,$
so that an exit does not occur at time $0$ even if $\Idx(X_{0})\neq\Idx(x).$

We will write $\mathbb{P}_x$ and $\mathbb{E}_x$ to be the respectively
the probability and expected value when the process starts from the
initial condition $x$. Similarly we will write $\mathbb{P}_\nu$ and
$\mathbb{E}_\nu$ to be   probability and expected value when the
initial condition is distributed as a the probability measure $\nu$.

We can now define the main object of our algorithm, the exit kernel
$Q$ defined on $A$ by 
\[
Q(x,\cdot)=\mathbb{P}_x(X_{\tau}\in\ccdot).
\]
 Note that $Q$ does not depend on $\kappa$, because the step where
 $\kappa$ is chosen is included in its definition.

Given any injection measure $\nu$ with $\nu(A_j)>0$ for all $j$, we define the associated weights $a_j(\nu)$
and stratified  injection measure $\nu_j$ by
\begin{align*}
  a_{j}=\nu(A_{j})\quad\text{and}\quad \nu_{j}(\cdot)=\frac{\nu(\;\cdot\;\cap A_{j})}{a_{j}}.
\end{align*}
 Note that $\nu=\sum_{j}a_{j}\nu_{j}$.
Next we define transition matrix $G$ by 
\[
G_{ij}=\mathbb{P}_{\nu_{j}}(X_\tau\in A_{i})=\nu_{j}Q(A_{i}).
\]

Finally, for each $j,$ define the exit measure from $A_{j}$
by 
\[
\xi_{j}(\cdot)=\mathbb{P}_{\nu_{j}}(X_{\tau}\in\cdot)=\nu_{j}Q(\cdot).
\]
Note that $G_{ij}=\xi_{j}(A_{i}).$

We will denote the above quantities associated with the associated
with the equilibrium injection measure $\nu^*$ by
$a_{j}^*,\nu_{j}^*,G^*$ and $\xi_{j}^*$. Notice that without loss of
generality, we can assume that $a_j^*>0$ for all $j$. Also, observe
that  $,G^*$ is the matrix which the index process in equilibrium
follows. In equilibrium, this is truly a Markov provess on the index
space since $\nu^*=\nu^* Q$.

In this framework, the basic version of our algorithm proceeds by\footnote{In practice, the algorithm will build a finite approximation of $\nu^{n+1}$
given $\nu^{n}$, by sampling a finite number of starting points
from $\nu^{n}$ and calculating an exit point from each.} $\nu^{n+1}=\nu^{n}Q$, and we will want to show that $\nu^{n}\rightarrow\nu^{*}$
as $n\rightarrow\infty.$

Analyzing the stratified processes will require understanding their
transition kernels. To this end, define, for each $j\in\{1,...,J\}$
and $x\in A_{j},$ the restricted kernel 
\[
\tilde{P}_{\kappa,j}(x,\cdot)=\frac{P_{\kappa}(x,\;\cdot\;\cap A_{j})}{P_{\kappa}(x,A_{j})}
\]
We will assume throughout that $P_{\kappa}(x,A_{j})>0$ if $x\in A_{j}.$
Define the corresponding un-normalized kernel 
\[
\hat{P}_{\kappa,j}(x,\cdot)=P_{\kappa}(x,\;\cdot\;\cap A_{j})
\]
and finally the un-normalized kernel restricted to leaving $A_{j}:$
\[
\bar{P}_{\kappa,j}(x,\cdot)=P_{\kappa}(x,\;\cdot\;\cap A_{j}^{c}).
\]
Denote the quasi-stationary distribution (QSD) of $P_{\kappa}$ on
$A_{j}$ by $\tilde{\nu}_{\kappa,j}.$ That is, $\tilde{\nu}_{\kappa,j}$
is the unique invariant distribution of $\tilde{P}_{\kappa,j}$. The
kernels $\hat{P}_{j}$, $\tilde{P}_{j}$, and $\bar{P}_{j}$ are defined
by integrating the associated kernel over $\kappa\sim\eta$, just as $P$
was defined from $P_{\kappa}$. 

The last notation we need to introduce has to do with the relationship
between injection measures and the invariant measure on the whole
space (and approximations of it), which is what we originally wanted
to sample from.

Given an injection measure $\nu_{j}$ on a stratum
$A_{j},$ define the corresponding occupation
measure, and it's normalization, by
\[
\mu_{j}(B)=\mathbb{E}_{\nu_{j}}\left[\sum_{k=0}^{\tau_{x}-1}{\bf
    1}_{(X_{k}\in B)}\right]
\qquad\text{and}\qquad
\pi_{j}(B)=\frac{1}{\mathbb{E}\tau_{\nu_{j}}}\mathbb{E}_{\nu_{j}}\left[\sum_{k=0}^{\tau_{\nu_{j}}-1}{\bf
    1}_{(X_{k}\in B)}\right],
\]
for $B \subset A$. Here $\mathbb{E}\tau_{\nu_{j}}$ means the expected value of the
exit time from $A_{j}$ for a particle started at $\nu_{j}.$ If we
also have weights $a_{j}$ for the strata, then we can define the
total occupation measure on $A$: 
\[
\mu=\sum_{j}a_{j}\mu_{j}\quad\text{and}\quad
\pi=\frac{1}{\mu(A)}\sum_{j}a_{j}\mu_{j}=\frac{1}{\mu(A)}\mu.
\]
(As before, $\pi_{j}^{*},\pi^{*},$ etc. are defined analogously). 

Throughout this article, $\parallel\cdot\parallel$ will mean $\parallel\cdot\parallel_{TV}$
unless otherwise stated. $a,a^{*}$will mean the vectors $(a_{1},...,a_{J}),$
$(a_{1}^{*},...,a_{J}^{*}),$ and similarly for other vectors and
matrices. We will also denote the weights given by a total injection
measure as $a(\nu)=(a_{1}(\nu),...,a_{J}(\nu))=(\nu(A_{1}),...,\nu(A_{J}))$.

\subsection{Statement of the Algorithm}

We are now ready to state how the injection measure method proceeds formally. 
The idea behind the algorithm is to start with some injection measures and weights, and on each iteration, 
calculate the exit measures given by them, then combine those into new exit measures and weights. 
The eigenvector version adds a step in which the starting weights are replaced by the 
eigenvector weights of the transition matrix given by the injection measures. 

For reference, we state the precise form of the injection measure
method in Algorithm 1

\begin{algorithm}
  \caption{The Injection Measure Method }
\begin{algorithmic}[1]
  \State  $N\gets \text{\# iterations}$
 \State $M\gets \text{\# Points per Strata} $
\State Initial weights $a_j^{0}$ and strata measures $\nu_{j}^{0}$
\Comment{Initial measure $\nu^{0}=\sum_{j}a_{j}^{0}\nu_{j}^{0}$}
\For{$n \gets\ 0\  \algorithmicto \ N-1$ }
 \For{$j \gets\ 1\  \algorithmicto \ J$ }
  \For{$i \gets\ 1\  \algorithmicto \ M$ }
	\State  $X_{j,i}^{n}(0) \gets \text{Random as }\nu_{j}^{n}$
        \State  $\kappa_{j,i}^n \gets  \text{Random         as }\eta$
        \While{$k < \tau$}
            \State $X_{j,i}^{n}(k+1)\gets \text{Random        as }  P_{\kappa_{j,i}^n}(X_{j,i}^{n}(k),\ccdot)$
            \EndWhile
           \State   $\mu_{j,i}^n  \gets \displaystyle\sum_{k=0}^{\tau-1} \delta_{X_{j,i}^n(k)}$
        \EndFor
       \State  $\mu_{j}^n \gets\frac1{M} \displaystyle\sum_{i=1}^M  \mu_{j,i}^n$
      \EndFor
      \For{$j,k \gets\ 1\  \algorithmicto \ J$ }
      \State $G_{kj}^{n} \gets \mid\{i:\Idx(X_{k,i}^n(\tau))=j\}\mid$
      \EndFor
      \State Normalize $G$ to be a probablity transition matrix

	\If{Basic~Version}
        \For{$j \gets\ 1\  \algorithmicto \ J$}
		\State $a_{j}^{n+1} \gets
                \displaystyle\sum_{k}a_{k}^{n}G_{kj}^{n}$\Comment{$a^{n+1} \gets a^n G^n$}

		\State $\nu_{j}^{n+1} \gets
                \frac{1}{a_{j}^{n+1}} \displaystyle\sum_{k}a_{k}^{n} \displaystyle\sum_{i:X_{k,i}^{n}(\tau)\in
                  A_{j}}\delta_{X_{k,i}^n(\tau)}$\Comment{New
                  injection measures}
                \EndFor

		\State $\nu^{n+1} \gets
                \displaystyle\sum_{k}a_{k}^{n+1}\nu_{k}^{n+1}\quad\text{and}\quad  \mu^{n+1} \gets \displaystyle\sum_{k}a_{k}^{n+1}\mu_{k}^{n+1}$
	\EndIf

	\If{Eigenvector~Version}
        	\State $z^n \gets \text{Normalize solution of }z^{n}G^{n}=z^{n}$\Comment{$z^{n}=(z^{n}_1,\dots,z^{n}_J)$}
		\State $\nu_{j}^{n+1}\gets \frac{1}{z_{j}^{n}}\sum_{k}z_{k}^{n}\displaystyle\sum_{i:X_{k,i}^{n}(\tau)\in
                  A_{j}}\delta_{X_{k,i}^n(\tau)}$\Comment{New
                  injection measures}
                	\State $\nu^{n+1}\gets
                        \sum_{j}z_{j}^{n}\nu_{j}^{n+1} \quad\text{and}\quad  \mu^{n+1} \gets \displaystyle\sum_{k}z_{k}^{n+1}\mu_{k}^{n+1}$
	\EndIf
\EndFor

\State Return : $\mu^{N}$ \Comment{Approximate of $\pi$}
\end{algorithmic}
\end{algorithm}

$\,$

\section{Main Results} \label{results}

\subsection{Assumptions}

Roughly, in order to prove our convergence theorems, we need to first assume that 
the process is well-behaved both within each stratum, and in terms of how it moves mass between the strata. 
Our strategy will then be to show that the exit process moves the weights closer to the true weights, at which point the behavior within each stratum lets us show that coupling can occur. 

The first two assumptions below can be thought of as non-degeneracy and
regularity conditions on $P$ and the strata. The first says, in effect,
that $P$ is well behaved and has a unique equilibrium, and that it
can move mass into each strata from at least one of the others. The
second is a Lyapunov-type condition we can use to control the growth
of the exit measure for a well-behaved injection measure, even when
some strata are unbounded.
\begin{assumption}[A0]
  \label{A0}  The Markov transition kernel $P$ has a unique invariant distribution $\pi$. 
  For all $\kappa$, $P_{\kappa}$ is irreducible and
has a unique invariant distribution $\pi_{\kappa}$, such that for
all $j$, 
\[
\intop_{A_{j}^{c}}P_{\kappa}(x,A_{j})\pi(dx)>0.
\]
and  $P(x,\cdot)=\int_{[0,1]}P_{\kappa}(x,\cdot)\eta(d\kappa)$.
\end{assumption}
\begin{assumption}[A1]
\label{A1} There exists a continuous function $V\colon A\rightarrow[0,\infty),$
a compact set $\mathcal{K}\subset A$, and $b,K>0,$ $\gamma\in(0,1)$
such that\footnote{If $V(x)\rightarrow\infty$ as $\mid x\mid\rightarrow\infty$, then
the condition is equivalent to requiring that $PV(x)\leq\gamma^{\prime}V(x)+b^{\prime}$
for some different $\gamma^{\prime}\in(0,1)$ and $b^{\prime}>0$. We will use
the slightly more general condition in Assumption~A1.} for all $x\in A$, $\kappa\in[0,1]$, 
\[
P_{\kappa}V(x)\leq\gamma V(x)+b{\bf 1}_{\mathcal{K}}
\]
and if $x\in\mathcal{K}$ and $\Idx(x)=j$, then 
\[
\mathbb{E}_{x}\tau_{j}\leq K.
\]
\end{assumption}
Next is our primary assumption on the system's ``microscopic'' behavior,
i.e. it's behavior inside a stratum. Intuitively, it states that any
particle starting from any point in a stratum has a chance of coupling
with the quasi-stationary distribution (QSD) of $P_{\kappa}$ in that
stratum before or at the time of exit. Therefore, it's exit distribution
will look like, with some probability, what it would have been had
the particle been injected according to  the QSD. 
\begin{assumption}[A2]
\label{A2} For all $j$ and $\kappa\in[0,1]$, there exists a
unique QSD $\tilde{\nu}_{\kappa,j}$ of $P_{\kappa}$ in $A_{j}$.
Let $\tilde{\nu}_{j}=\int_{[0,1]}\tilde{\nu}_{\kappa,j}\eta(d\kappa)$,
where $\eta$ is the distribution $\kappa$ is drawn from at the start
of a trajectory. Then there exists a constant $0<c<1$ such that,
for any $j$ and $\nu_{j}^{0},$ 
\[
\xi_{j}^{0}\geq c\cdot\tilde{\xi}_{j}
\]
where $\tilde{\xi}_{j}$ is the exit measure from $A_{j}$ started
from $\tilde{\nu}_{j}$, via the kernel $Q$. That is,  $\tilde{\xi}_{j}=\tilde{\nu}_{j}Q$.
\end{assumption}
Next, we need an assumption on the ``macroscopic'' behavior, i.e.
how mass moves between strata. (Assumption~A0 ensures some
movement between strata, but not does not ensure that it is global or
give quantitative  information.)  We have three versions of this assumption,
and a convergence theorem that holds under each.c The first form says
that after enough exits, a particle has a probability, bounded from
below, of being in any of the strata. i.e. the exit process explores
the space after enough time. The second choice of assumption says
that the transition matrix in equilibrium, $G^{*}$, is geometrically
ergodic, with mixing rate $\lambda$. This will allow us the prove
a convergence result for the basic version with a more precise estimate
of the convergence rate, combining the microscopic rate $c$ and macroscopic
rate $\lambda$.

Note that the above assumptions could fail due to the index of a particle
on the $n$-th exit being periodic. For example, if there are only two strata, then 
the index will keep flipping between one stratum and the other between exits. If this is the only 
problem keeping the assumptions from holding, it can be remedied by introducing laziness to the 
exit kernel. That is, $Q$ can be replaced by $pQ+(1-p)I$, where $I$ is the identity kernel, and
the periodicity will be removed. 

The last form of the macroscopic assumption pertains to the eigenvector
version, and is borrowed from the perturbation results in \cite{thiede_koten_weare_2015}.
Instead of giving the mixing rate of $G^{*}$, we assume that the
invariant distribution vector $G^{*}$ has bounded sensitivity to
perturbations in its entries. We will then prove a convergence theorem
for the eigenvector version, this time combining $c$ and the sensitivity
constants. We state all three macroscopic assumptions below. 
\begin{assumption}[B1]
There exists $u>0$ and $m\geq1$ such that, for any $j,$ if
$X_{0}\sim\tilde{\nu}_{j}$ and $X_{n}$ is the position of the $n$-th
exit from one stratum to another, then 
\[
\Law(\Idx(X_{n}))\geq ua^{*}
\]
for all $n\geq m.$ 
\end{assumption}
\begin{assumption}[B2]
There exists $m\geq1$ and $\lambda^{*}\in(0,1)$ such that,
for $n\geq m$ and any two probability vectors $a,b$ on $\{1,...,J\},$
\[
\parallel a(G^{*})^{n+1}-b(G^{*})^{n+1}\parallel\leq\lambda\parallel a(G^{*})^{n}-b(G^{*})^{n}\parallel.
\]
\end{assumption}
\begin{assumption}[B3]
 There exists $\theta_{ik}>0$ for $i,k\in\{1,...,J\}$,
depending only on $G^{*}$ such that, if $G$ is an irreducible transition
matrix on $\{1,...,J\}$ such that, if $G,G^{*}\geq c\tilde{G}$\footnote{Where $\tilde{G}_{ij}=\tilde{\xi}_{j}(A_{i})$ is the transition matrix
between strata when the injection measures are the QSD's.} and $z,z^{*}$ are the invariant measures of $G,G^{*},$ resp., then
for all $j,$ 
\[
  \sup_{j}\left(\log(z_{j})-\log(z_{j}^{*})\right)\leq\sum_{i\neq k}\theta_{i,k}\mid G_{ik}-G_{ik}^{*}\mid
\]
\end{assumption}

\subsection{Main Theorems}

We are now ready to state our main theorems, as well as give a brief
strategy for proving them. The detailed proofs will be reserved for
section 5.

The first theorem states in effect that the equilibrium injection
measure exists uniquely, and that the corresponding occupation measure
is the original measure we wanted to sample. Our strategy for proving
this involves manipulating the sub-stochastic kernels $\hat{P}_{\kappa,j}$
and $\bar{P}_{\kappa,j}$, and how they relate the injection and occupation
measures. It is essentially a Poisson equation argument at heart:
$\mu_{j}$ solves a Poisson-like equation in terms of $\nu_{j}$,
and combining these gives an equation for $\mu$ that reduces
to $\mu(I-P)=0$ in the equilibrium case.

\begin{theorem}\label{Theorem 1} Suppose that (A0) and (A1) hold.
Then there exists a unique probability measure on $A,$ $\nu^{*}$
such that $\nu^{*}V<\infty$ and $\nu^{*}Q=\nu^{*}.$ Furthermore,for
the coresponding occupation measure, $\mu^{*}=\mu^{*}P,$ where $P(x,\cdot)=\int_{[0,1]}P_{\kappa}(x,\cdot)\eta(d\kappa)$.
\end{theorem}

Our second main theorem establishes convergence of the basic version
after enough steps are performed. We use a fairly straightforward
coupling argument to prove it. The idea is that assumption (B1) states
that two particles following the exit process $Q$ can eventually
be in the same stratum, at which point (A1) says that they have a
chance of being coupled after the next iteration of $Q$. This establishes
geometric ergodicity of $Q$.

\begin{theorem} \label{Theorem 2} Under (A1) and (B1), for any $n\geq m$,
$Q^{n+1}$ is a global contraction on probability measures on $A$,
with contraction constant $1-c^{2}u$. In particular, $\parallel\nu^{k(m+1)}-\nu^{*}\parallel\leq(1-c^{2}u)^{k}$
for $k\geq1$. \end{theorem}

Next, we have a theorem giving a more precise rate of geometric convergence
in the long term. To prove it, we will need to construct a new metric
which balances the total variation distance between $\nu^{n}$ and
the true $\nu^{*}$, and the difference between their strata weights.
The idea is that, even if no coupling is possible on a given step,
the weights will still get closer to the truth, which will allow coupling
with $\nu^{*}$ to occur at some future time. Since the metric we
show contraction in is equivalent to TV-distance, we get the desired
result.

\begin{theorem} \label{Theorem 3} Suppose that (A0)-(A2) and (B1)-(B2)
hold. For any $\nu^{0}$ such that $\nu^{0}V<\infty$, there exists
$q_{n},$ $\forall n\geq0,$ such that $\parallel\nu^{n+1}-\nu^{*}\parallel\leq q_{n}\parallel\nu^{n}-\nu^{*}\parallel$
and as $n\rightarrow\infty$, 
\begin{equation}
q_{n}\rightarrow q:=\inf_{\beta\in(0,1)}\left(\inf_{0<\alpha<\frac{\beta c}{S(1-c)}}\max\left(1-\beta c+\alpha S(1-c),\frac{1+(1-\beta)\alpha\lambda}{1+(1-\beta)\alpha}\right)\right)\label{eq:1}
\end{equation}
where $S=\frac{1}{1-\lambda}$. \end{theorem}

The expression we get for the limiting rate is complex, but has a
meaningful interpretation. Note that it takes the form of a maximum
of a rate in terms of $c$ and one in terms of $\lambda$. Therefore,
it suggests that, of the macro- and microscopic rates, whichever is
slower acts as a sort of bottleneck. It can also be shown straightforwardly
that if $1-c=0$, then the limiting rate is $\lambda$, and if $\lambda=0$,
then the limiting rate is order $\frac{1}{2}$ in $1-c$. So we can
conjecture that, to first order, it behaves like $\max(\sqrt{1-c},\lambda)$.

Our final theorem establishes geometric convergence of the eigenvector
version, for starting guess $\nu^{0}$ sufficiently close to $\nu^{*}$.
The final rate involves the microscopic coupling parameter $c$ and
a constant $r^{\infty}$ relating to how sensitive $G^{*}$'s eigenvector
is to perturbations. In this way, it is like an analogue to the previous
theorem, with $\lambda$ replaced by $r^{\infty}$. It also does not
have the bottleneck form, suggesting that if $\lambda$ is a slow
rate, then the eigenvector version may be faster than the basic version,
which is to be expected.

\begin{theorem}\label{Theorem 4} Let $\nu^{0},$ $\nu^{1},...$
be the total injection measures for a run of the eigenvector version.
Suppose that (A0)-(A2) and (B3) hold. Let 
\[
r^{\infty}=2(1-c)\sup_{j}(a_{j}^{*})\sup_{i,k}\left(\frac{e^{\theta_{ik}}-1}{a_{i}^{*}}\right)
\]
\[
E=exp\left(2(1-c)\sum_{i}\sup_{k\neq i}(\theta_{ik})\right)
\]
Also suppose that 
\[
\parallel\nu^{0}-\nu^{*}\parallel<\frac{1+r^{\infty}E}{2(1-c)\sup_{i\neq k}{\frac{\theta_{ik}}{a_{i}^{*}}}}\left(\frac{1}{q}-1\right),
\]
where $q=1-\frac{c}{r^{\infty}E+1}$. Then there exists $p_{n},$
$\forall n\geq0,$ such that $\parallel\nu^{n+1}-\nu^{*}\parallel\leq p_{n}\parallel\nu^{n}-\nu^{*}\parallel$
and
\[
p^{n}\rightarrow1-\frac{c}{r^{\infty}+1}.
\]
\end{theorem}

\section{Numerical Simulations} \label{numerics}

We now turn to analyzing the results of numerical simulations of our
algorithms. We focus mainly on the eigenvector version here, as in
all our simulations its convergence is strictly faster than the basic
version. but does not immediately suggest how to prove this. We include
the code used to generate our results in the following git repository:
https://gitlab.com/gabeearle/julia-code-stratified-2021

We test our method on the two-dimensional Maier-Stein system, as outlined
in \cite{maier_stein_1996}, with parameters used in \cite{heymann_2007}.
This is a relatively simple low-dimensional system, but one which
displays the non-reversibility which our methods are suited for. Furthermore,
with the chosen parameters, the system displays both increased non-gradient
effects and a double-well type invariant distribution.

For reference, we state the version of the Maier-Stein system we use
here. The dynamics evolve according to the following SDE: 
\[
du=(u-u^{3}-\beta uv^{2})dt+\sqrt{\epsilon}dW_{t}^{1}
\]
\[
dv=-(1+u^{2})vdt+\sqrt{\epsilon}dW_{t}^{2}
\]
where $\beta,\epsilon>0$ and $W^{1},W^{2}$ are independent Brownian
motions. For $\beta\neq1$, this system is non-reversible. For $\beta>4$,
it displays additional non-gradient behavior and unusual minimum action
paths, as detailed in \cite{heymann_2007}. We use the parameters
$\beta=10$ and $\epsilon=0.01$. Fig. 1 Shows an approximate 2D histogram
plot of the invariant distribution of this system, computed via a
discretization of the dynamics.

In Fig. 2, we show several ways to choose the strata. One is a subdivision
of the space into 3 vertical ellipses, another into 5 smaller ones.
The next is a version in which the 5 strata are rotated, so that they
do not perfectly line up with the axes of the system. Finally, we
also study the case where the space is divided into 6 circular strata,
with varying sizes and more than two strata overlapping at once. 

\label{Fig. 1} 
\begin{figure}
\begin{centering}
\includegraphics[scale=0.4]{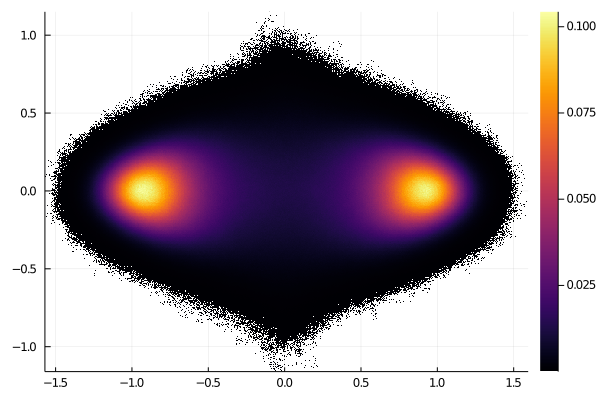}
\par\end{centering}
\caption{A 2D histogram of the invariant measure for the Maier-Stein system.
This approximation of the true density was generated by an (un-stratified)
Runge-Kutta Method, with noise, run for $10^{8}$ steps, and averaged
over two runs. The $u$-marginal density of this distribution will
be used as our benchmark for the true density when calculating the
error of the stratified method in $u$.}
\end{figure}

\label{Fig. 2} 
\begin{figure}
\begin{centering}
\includegraphics[scale=0.3]{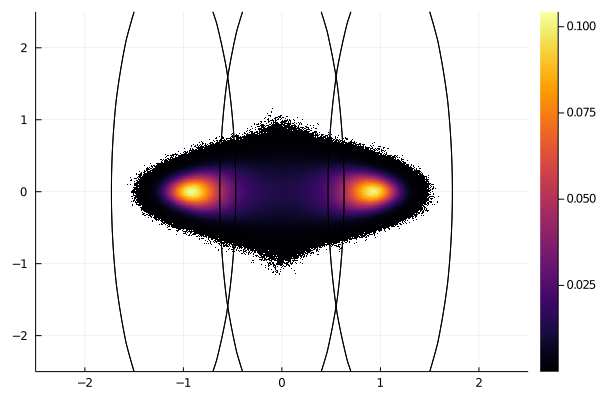}\includegraphics[scale=0.3]{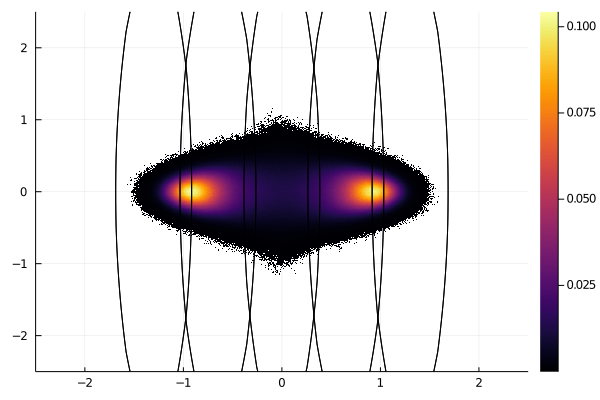}
\par\end{centering}
\begin{centering}
\includegraphics[scale=0.3]{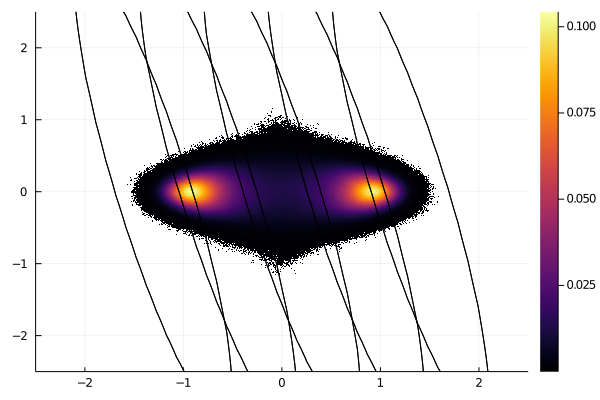}\includegraphics[scale=0.3]{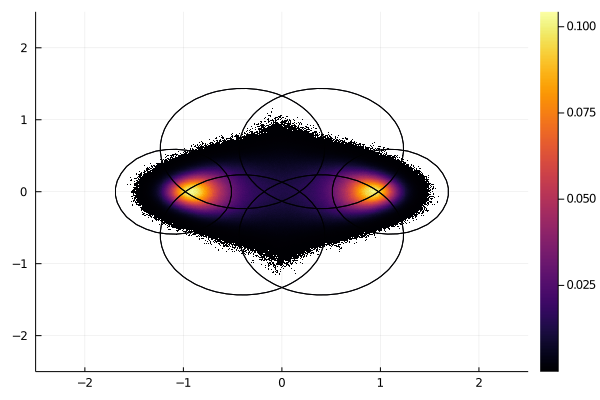}
\par\end{centering}
\caption{Four choices of elliptical strata for the Maier-Stein system. We use
setups with 3 and 5 vertically oriented strata, tilted strata, and
6 circular strata which cover the space more tightly}
\end{figure}

We can begin by comparing the 1-dimensional projections of both estimates
to get a sense of how well the algorithm approximates the double-well
structure of the system. Fig. 3 shows the un-stratified and stratified
(with 3000 exits) estimates of the $u$-marginal after 30 iterations.
As we can see, they agree quite well, demonstrating the accuracy of
our method.

\label{Fig. 3} 
\begin{figure}
\begin{centering}
\includegraphics[scale=0.4]{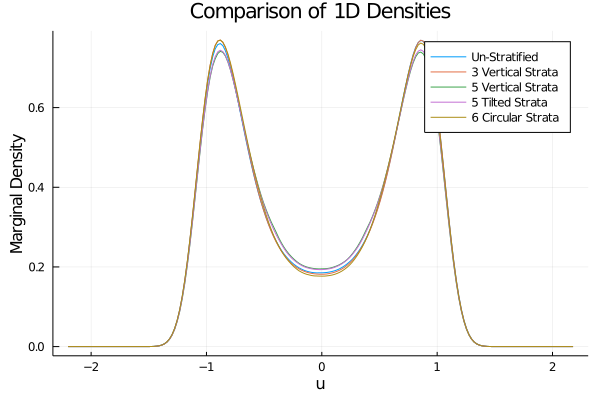}
\par\end{centering}
\caption{1-dimensional histograms for the true $u$-marginal density
  (again from a high-resolution un-stratified run) and the
approximations of it generated by the algorithm with the strata setups
above. Note that, with some deviation, each setup produces quite an
accurate approximation of the marginal density.}
\end{figure}

In Fig. 4 we plot, for each of the strata setups, the final total
injection measure, i.e. our approximation of $\nu^{n}$, where
$n=30$ is our final iteration number.Note that in each case, in injection measure has converged to a density roughly along the boundaries of the ellipses, illustrating that the algorithm can be used to estimate this density as well as the invariant distribution of the system. 

\label{Fig. 4} 
\begin{figure}
\begin{centering}
\includegraphics[scale=0.3]{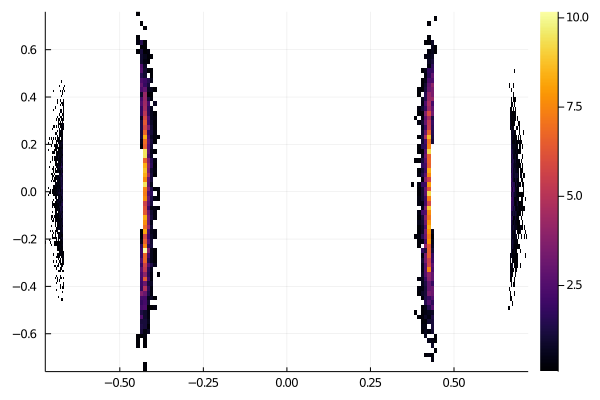}\includegraphics[scale=0.3]{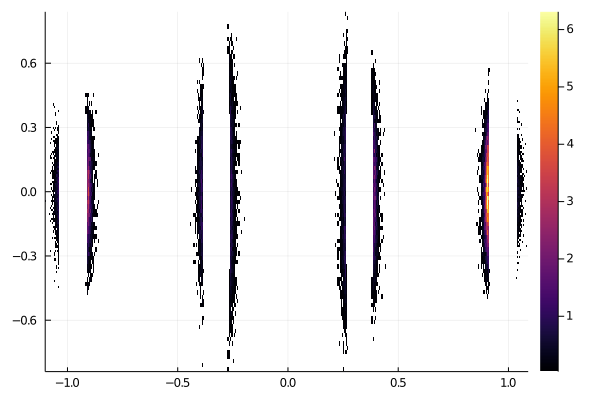}
\par\end{centering}
\begin{centering}
\includegraphics[scale=0.3]{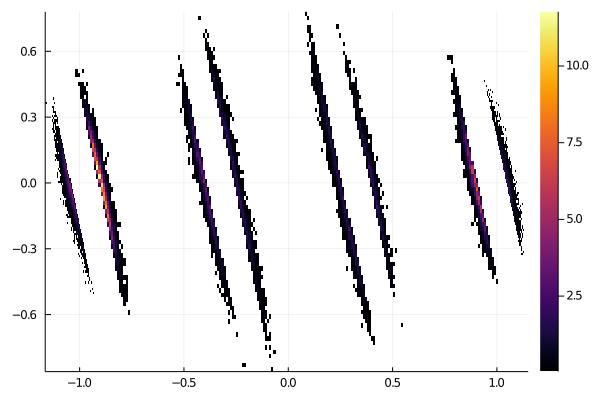}\includegraphics[scale=0.3]{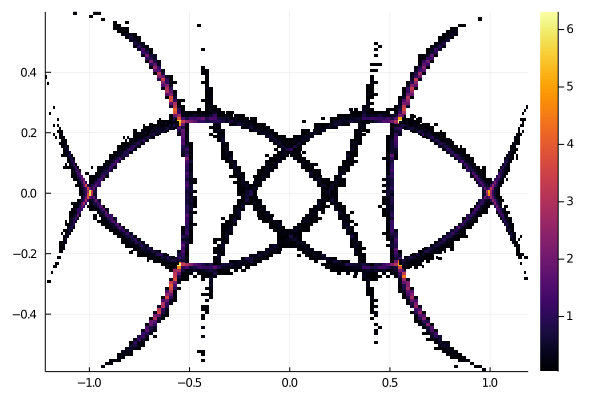}
\par\end{centering}
\caption{2D Histogram of the density of the total injection measure, for the
strata setups in Fig. 2. The densities are calculated from the points
generated on the last iteration of a 30-iteration run of the algorithm,
with 3000 exits per iteration. Because we are using strata with hard
boundaries, the injection measures lie on a portion of the boundaries
of each ellipse.}
\end{figure}

We can also examine how much the number of exits per iteration affects
performance. Fig. 5 shows a log-log plot of the error in the $u$-marginal
of the invariant distribution vs. time, for a range of exits the algorithm
was run with and for each setup. The error is approximately computed
by taking the TV-distance between the histogram from the un-stratified
run shown in Fig. 3, and the histogram computed from the stratified
run, averaged up to up to iteration $n$. The time spent is computed
as the total number of steps of the SDE taken in all strata up to
iteration $n$. The results, are shown in Fig. 5.

We can see that, in the 3-strata setting, increasing the number of
exits (hence the number of sample from each exit measure) per iteration
does not seem to affect the accuracy per computational time greatly.
The runs of the algorithm with more exits do in fact decrease the
error more quickly on each iteration, but the iterations also take
longer, and in the end the effects balance out, and all three runs
have similar trend lines of error vs. time. For both 5-strata setups,
this behavior changes dramatically. The version of the algorithm run
with 3000 exits still takes longer to perform, but also sees error
decrease much more rapidly. The reason for such an increase in accuracy
is not immediately clear, and could be a subject of future work. The
interesting phenomenon is that there appears to be some threshold,
in terms of the work put in on each iteration, at which the compounding
error becomes much less of an issue. 

\label{Fig. 5} 
\begin{figure}
\begin{centering}
\includegraphics[scale=0.3]{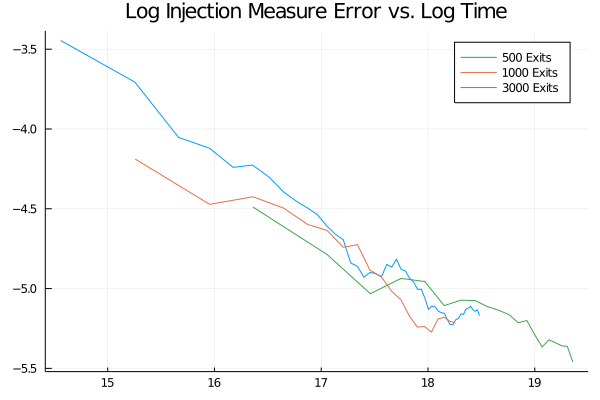}\includegraphics[scale=0.3]{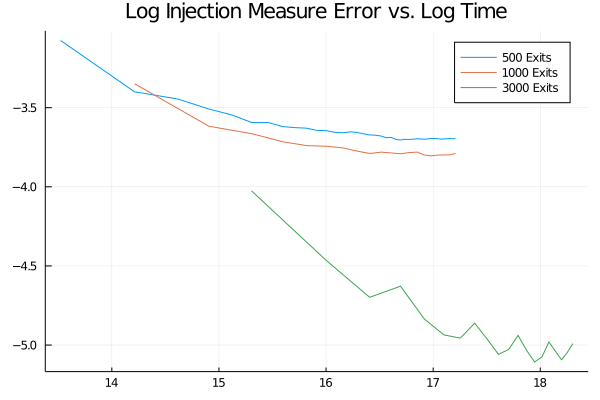}
\par\end{centering}
\begin{centering}
\includegraphics[scale=0.3]{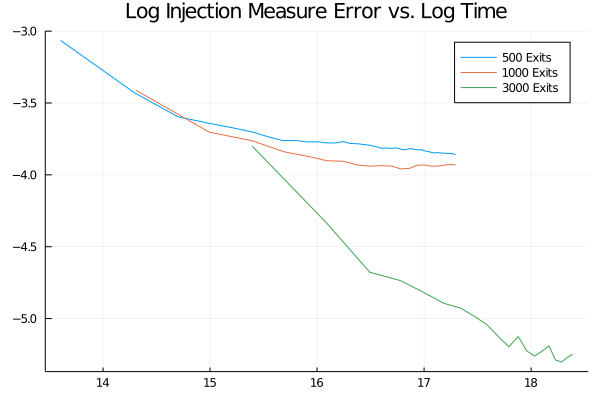}\includegraphics[scale=0.3]{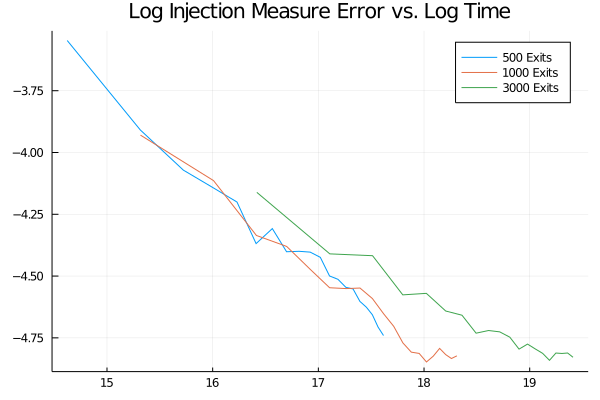}
\par\end{centering}
\caption{Log-Log total variation error in $u$ vs time, for each of the strata
setups. The error values are averaged over 20 full runs of the algorithm
for each setup. Computational time is measured as the total number
of steps of time $h$ that each particle has taken, up to the current
iteration. Note that, for the two 5-strata setups, a dramatic increase
in accuracy occurs somewhere between 1000 and 3000 exits run per iteration.}
\end{figure}

A possible hint as to why the accuracy threshold occurs can be found
by examining how much the injection measure changes for each strata
setup and each choice of the number of exits. We measure fluctuations
in the total injection measure similarly to the error in the invariant
measure: by taking the TV-distance between a histogram of the approximate
$\nu$ marginal in the $v$-direction on one iteration, and
the same histogram on the next iteration. The results are shown, for
all the setups, in Fig. 6. First, we can see that the total variation
change in the injection measure from one iteration to the next remains
quite small for each case. However, corresponding to the large drops
in error fore the 5-strata setups are cases where, when going from
1000 to 3000 exits, the fluctuation not only decreases but seems to
behave periodically. What this means for the algorithm's deeper behavior
is not immediately clear, but it could be related to the accuracy
threshold observed in Fig. 5. 

\label{Fig. 6}

\begin{figure}
\begin{centering}
\includegraphics[scale=0.3]{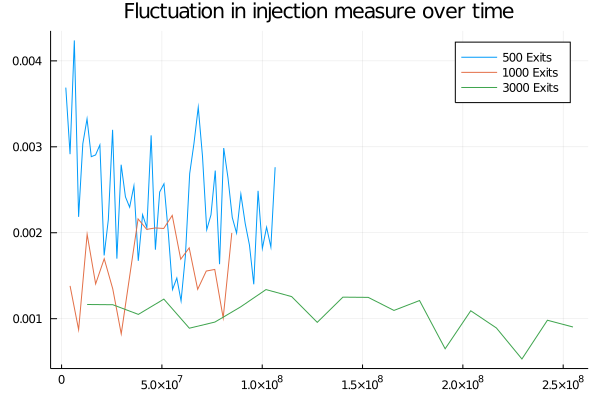}\includegraphics[scale=0.3]{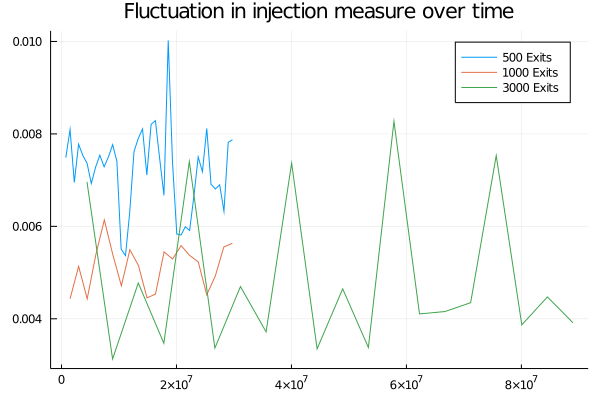}
\par\end{centering}
\begin{centering}
\includegraphics[scale=0.3]{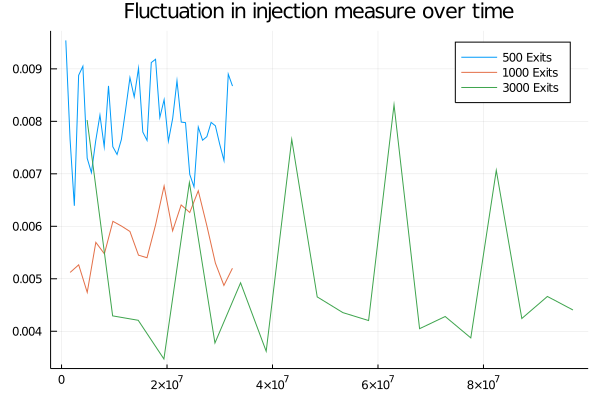}\includegraphics[scale=0.3]{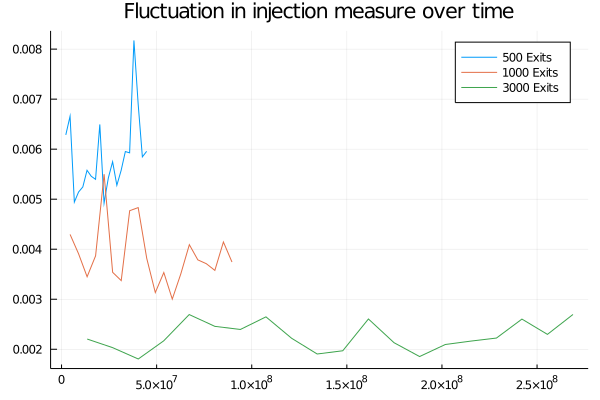}
\par\end{centering}
\caption{Fluctuation of the injection measure vs. time for the strata setups.
The fluctuation is computed as the TV distance between the $v$-marginal
of the total injection measure on one iteration and on the next. Computational
time is computed as in Fig. 6. Corresponding to the jump in accuracy
in the 5-strata setups is a lower and more regular fluctuation in
the injection measure.}
\end{figure}

In Fig. 7, we plot the number of points in the injection measure
of each stratum, for each setup and some choices of the number of
exits per iteration. Unsurprisingly, we see in the 5-strata case that
having more exits results in the more points in each injection measure.
More specifically, we can see a large difference in how many points
are in the ``smallest'' of the strata. Perhaps the accuracy threshold
has some connection to there being enough exits that the least common
strata still gets enough points to approximate $\nu_{j}^{n}$ reasonably
well. In that case, an alternate version of the algorithm could be
devised, in which more exits are sampled from the strata neighboring
the ``smallest'' ones, as that is where it is most important to
get samples. This would have the advantage of avoiding unnecessary
computational work in the strata which naturally get many points injected
already. We can already see from Fig. 7 that one of the main advantages
of stratified sampling holds here: that the different strata have
different sizes, in terms of how likely particles are to enter them,
but we put in similar amounts of work in them, by always drawing the
same number of points from each injection measure.

\label{Fig. 7}

\begin{figure}
\begin{centering}
\includegraphics[scale=0.3]{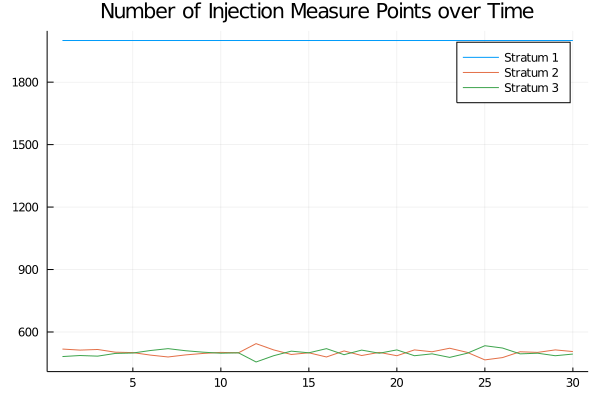}\includegraphics[scale=0.3]{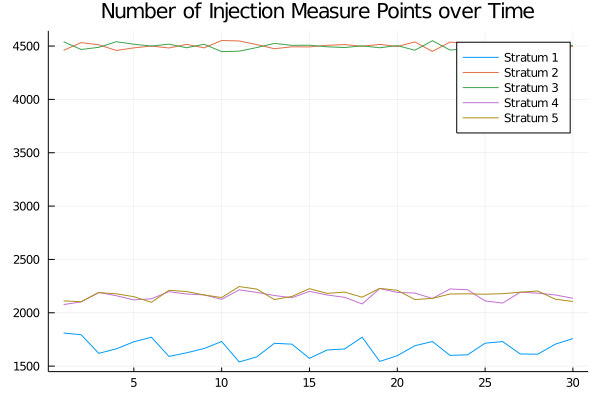}
\par\end{centering}
\begin{centering}
\includegraphics[scale=0.3]{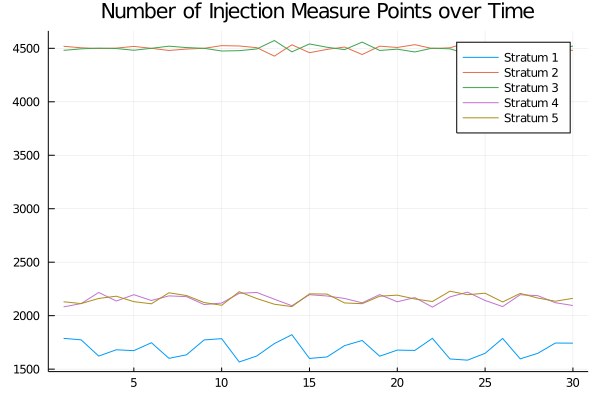}\includegraphics[scale=0.3]{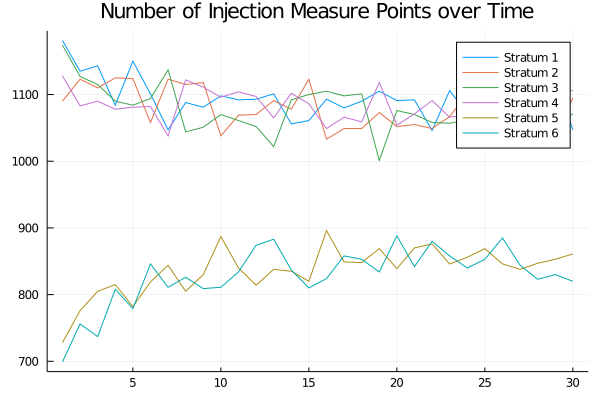}
\par\end{centering}
\caption{Injection measures sizes vs. iteration number for the strata setups.
The injection measure size is the total number of points that exited
into a given stratum on the current iteration. The figure illustrated
that the size of each measure converges to some value and then fluctuates
around it. Note that some strata have significantly more points injected
than others, but the same number of injections are drawn and run to
exit on an iteration for each, illustrating that similar amounts of
work are done in each strata despite their different sizes.}
\end{figure}

We can derive some insight into how the error behaves by decomposing
into the error in the weights or in each of the strata. In doing so,
we can see how the error in the weights relates to the overall error,
and in each stratum. In Fig. 8, We first plot the ``occupation weights''
of the strata for the 5 vertical strata setup, that is, their weight
as a fraction of the total mass in the invariant measure. We then
plot the error in the occupation measure for each stratum. For this,
we use as a benchmark the occupation measure computed from an un-stratified
run, where the measure in each stratum is formed from all points visited
in that stratum over the run. The error is then taken to be the total
variation distance between a histogram of the stratified run's occupation
measure and this benchmark. Finally, in Fig. 8 we also show in error
in the injection weights over time and the overall total variation
error over time.

Comparing the results for stratified runs with different numbers of
exits per iteration, we see that the overall error behaves almost
identically to the error in strata 2 and 3. This is not surprising,
as those are the strata with by far the highest occupation weight.
Since they dominate the mass of the invariant measure, most of the
total error will be composed of the error in them. We can see that
the error in the dominant strata sees a dramatic drop somewhere between
1000 and 3000 exits, suggesting that the same behavior we saw before
in the total error is caused by the drop occurring in strata 2 and
3.

Another observation to note is that the error in the weights also
behaves very similarly to the error in the dominant strata. A sharp
increase in accuracy in those occupation measures when going from
1000 to 3000 exits corresponds to a similar rise in accuracy in the
injection weights. Recall that the injection weights at iteration
$n$ are given by the eigenvector of the transition matrix $G^{n}$,
which itself is a function of the injection measures and the exit
kernel. The occupation and injection measures are also directly related
to each other. The behavior in the weights being similar to that of
strata 2 and 3 suggests that, for this choice of strata, the weights
are especially sensitive to perturbations in those strata. Therefore,
it seems that the weights are especially sensitive to perturbations
in the corresponding entries in the $G$ matrix.

Hence, in this case, the strata which are dominate the invariant measure's
mass and those to which the eigenvector weights are most sensitive
are the same.

\label{Fig. 8} 
\begin{figure}
\begin{centering}
\includegraphics[scale=0.2]{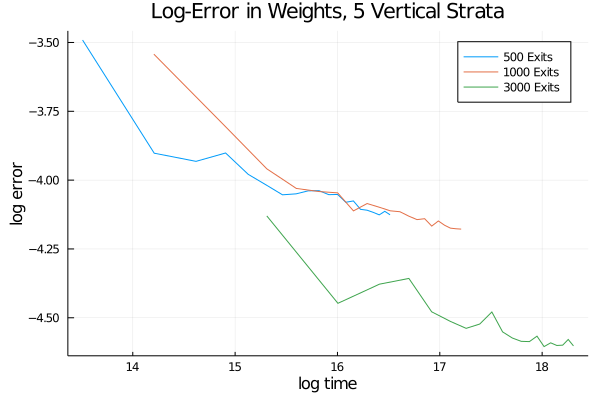}\includegraphics[scale=0.2]{5_25_MS_avged_exits_flux_5_reg_large_strata_30_runs_log_error_plot}\includegraphics[scale=0.2]{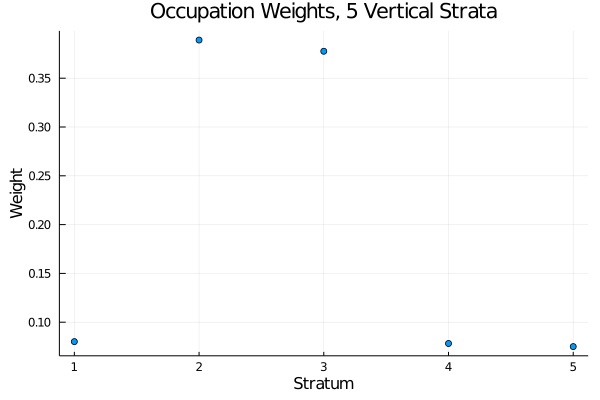}
\par\end{centering}
\begin{centering}
\includegraphics[scale=0.2]{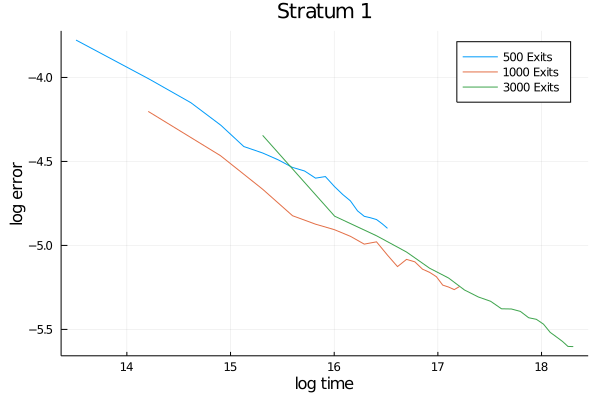}\includegraphics[scale=0.2]{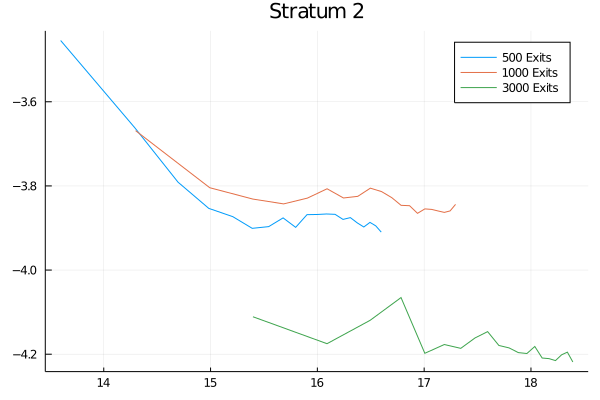}\includegraphics[scale=0.2]{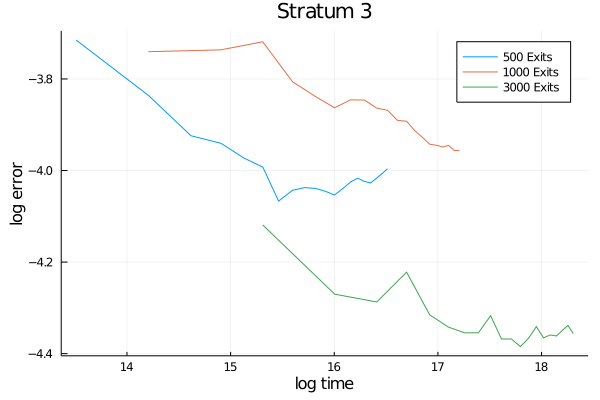}
\par\end{centering}
\begin{centering}
\includegraphics[scale=0.2]{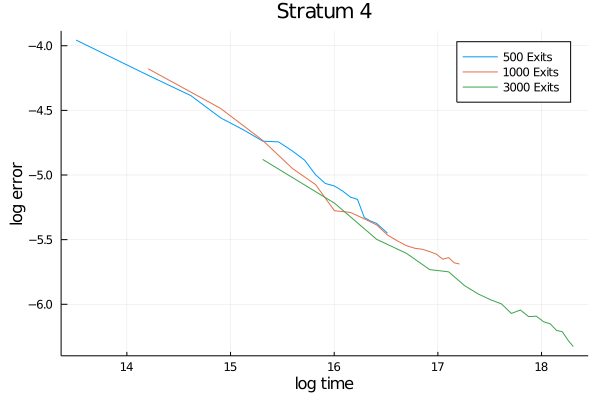}\includegraphics[scale=0.2]{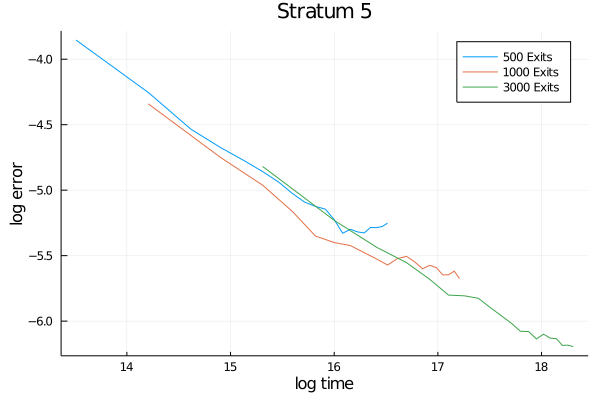}
\par\end{centering}
\caption{Occupation weights, individual strata errors, weights error and total
error, for 5 vertical strata setup. The error in individual occupation
measures is calculated as the TV distance between the $u$-histogram
generated by an un-stratified run and the $u$-histogram for the occupation
measure from the stratified algorithm, with the same bins as for the
overall error. The total error behaves similarly to the error in the
largest-weight strata, as does the weight error.}
\end{figure}

On the other hand, for a different choice of strata, we see slightly
different behavior. In Fig. 9, we plot the occupation weights, error
in occupation measures and weights, and total error in the same way
as before, but for the 6 circular strata setup. As before, we can
see that the total error behaves similarly to that of the highest
weight strata, this time strata 5 and 6. However, the injection weight
error behaves more like the occupation error in the previous 4 strata,
suggesting that the transition matrix $G$ for this setup is more
sensitive to the entries corresponding to transitions from those strata.
So as before, some strata are more important to the overall error
and some are more important than the weights, but this time the two
do not line up. Strata with small occupation weight have a large influence
on the accuracy of the eigenvector weights. This is important to
note, as it illustrates one of the main motivations for our method:
that strata which are visited rarely may still have outsize influence
on something of interest. 

\label{Fig. 9} 
\begin{figure}
\begin{centering}
\includegraphics[scale=0.2]{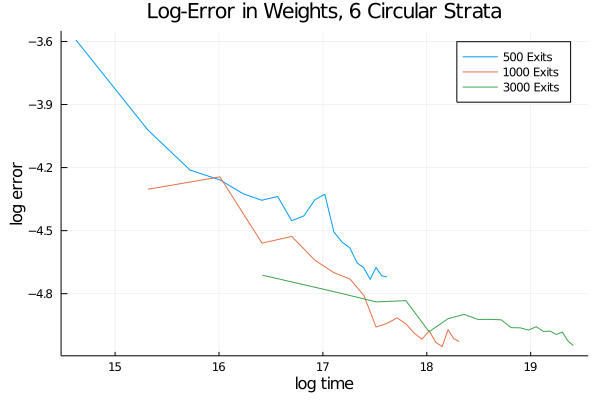}\includegraphics[scale=0.2]{6_27_MS_avged_exits_flux_6_circ_strata_30_runs_log_error_plot}\includegraphics[scale=0.2]{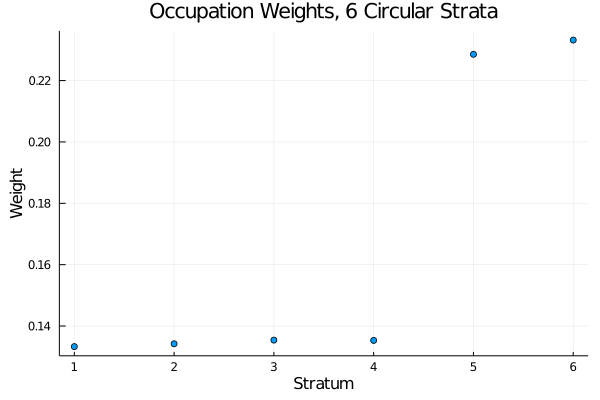}
\par\end{centering}
\begin{centering}
\includegraphics[scale=0.2]{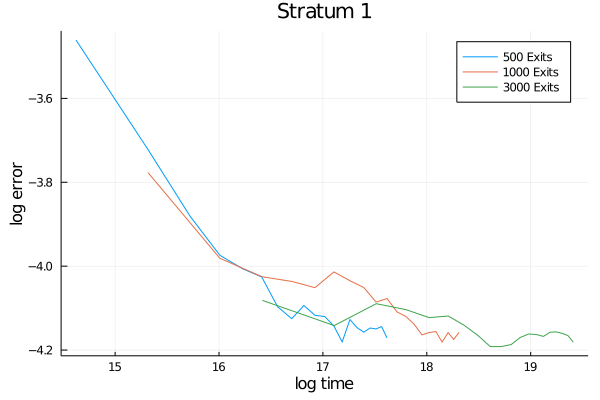}\includegraphics[scale=0.2]{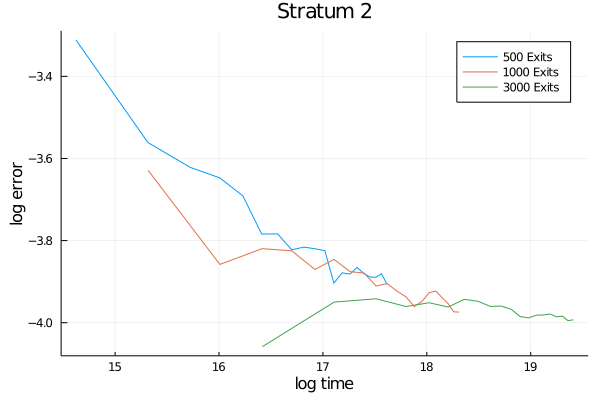}\includegraphics[scale=0.2]{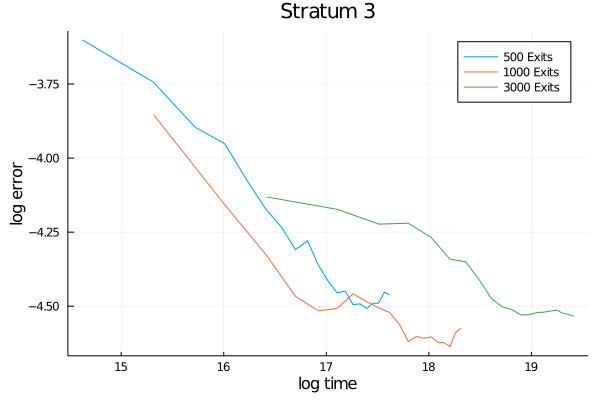}
\par\end{centering}
\begin{centering}
\includegraphics[scale=0.2]{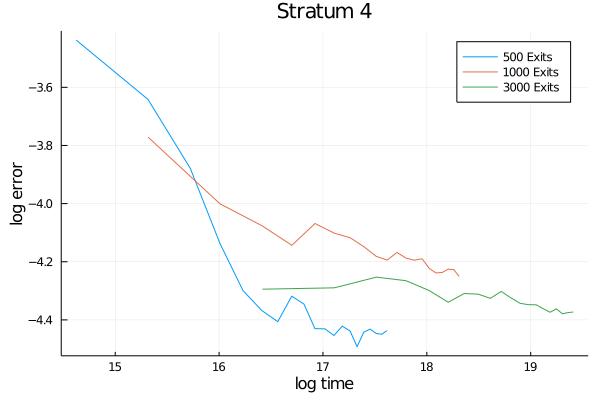}\includegraphics[scale=0.2]{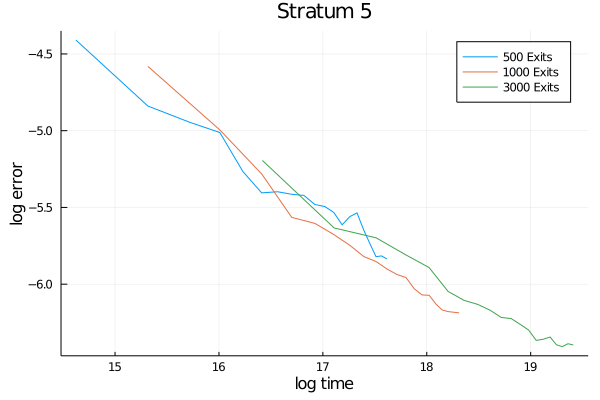}\includegraphics[scale=0.2]{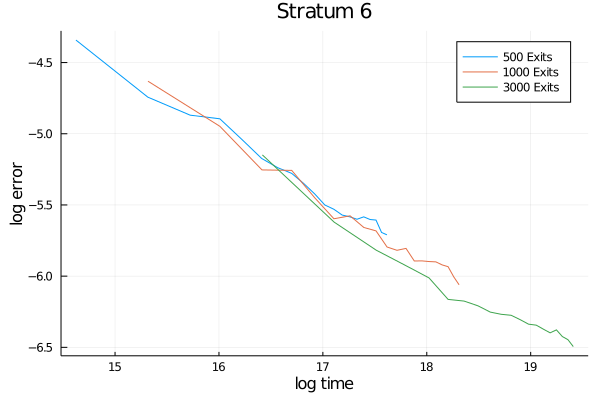}
\par\end{centering}
\caption{Occupation weights, individual strata errors, weights error and total
error, for the 6 circular strata setup. Errors are calculated the
same way as in Fig. 10. In this case, the total error still behaves
similarly to the error in the largest strata, but the weight error
is more similar to that of the smaller strata, showing that low-weight
strata can have significant importance.}
\end{figure}

\section{Proofs} \label{proofs}

\subsection{Existence and Uniqueness of the Equilibrium Injection Measure $\nu^{*}$}

The first result we must prove is that the fixed point of $Q$ exists,
and has the appropriate distribution. We will show this by manipulating
the sub-stochastic kernels derived from $P_{\kappa}$, and using the
ways they relate the injection, exit and occupation measures. First,
however, we must establish some regularity. To this end, our first
lemma states that, if a given injection measure has controlled growth
with respect to the Lyapunov function in assumption (A1), then the
corresponding exit and occupation measures exist, and also have controlled
growth. 

\begin{lemma} \label{Lemma 1} Suppose that (A1) holds. Let $\nu_{j}$
be a probability measure on $A_{j}$, with exit time, exit measure
and occupation measure $\tau_{j}$, $\xi_{j}$ and $\mu_{j}$, such
that $\nu_{j}V<\infty$. Then $\mathbb{E}_{\nu_{j}}\tau_{j}<\infty$,
$\xi_{j}$ and $\mu_{j}$are well defined, and $\xi_{j}V<\infty$,
$\mu_{j}V<\infty$. Furthermore, if $\pi$ is the invariant distribution
of $P_{\kappa}$ for any $\kappa\in[0,1]$, then $\pi V<\infty$ and
$\pi\mid_{A_{j}}(I-\hat{P}_{j})<\infty$. \end{lemma}

\begin{proof}  Fix $\kappa\in[0,1]$. Let $B\subset\mathcal{K}^{c}$
where $\mathcal{K}$ is defined in A1, and let $\tau_{B}=\inf\{n:X_{n}\notin B\}$
Where $X_{0}=x$ and $X_{n}$ follows $P_{\kappa}$. Define 
\[
V_{B}^{*}=\inf\{V(x):x\in B\}
\]
Following ideas from \cite{MeynTweedie93,MattinglyStuartHigham02}, we define the following process: 
\[
M_{n}=\frac{V(X_{n})}{V_{B}^{*}\gamma^{n}}-b\sum_{k=0}^{n-1}\frac{{\bf 1}_{(X_{k}\in\mathcal{K})}}{\gamma^{k+1}}
\]
 with the filtration $\mathcal{F}_{n}=\sigma(X_{n}:k\leq n)$. Then,
using the condition of (A1) to bound $\mathbb{E}\left[V(X_{n+1})\mid\mathcal{F}_{n}\right]=PV(X_{n})$,
we can simplify 
\[
\mathbb{E}\left[M_{n+1}\mid\mathcal{F}_{n}\right]=\frac{\mathbb{E}\left[V(X_{n+1})\mid\mathcal{F}_{n}\right]}{V_{B}^{*}\gamma^{n+1}}-b\mathbb{E}\left[\sum_{k=0}^{n-1}\frac{{\bf 1}_{(X_{k}\in\mathcal{K})}}{\gamma^{k+1}}\mid\mathcal{F}_{n}\right]-b\mathbb{E}\left[\frac{{\bf 1}_{(X_{n}\in\mathcal{K})}}{\gamma^{n+1}}\mid\mathcal{F}_{n}\right]
\]
 to show that 
\[
\mathbb{E}\left[M_{n+1}\mid\mathcal{F}_{n}\right]\leq\frac{V(X_{n})}{V_{B}^{*}\gamma^{n}}-b\sum_{k=0}^{n-1}\frac{{\bf 1}_{(X_{k}\in\mathcal{K})}}{\gamma^{k+1}}=M_{n}.
\]
 Furthermore, if $\mathbb{E}\mid M_{n}\mid<\infty$, then $\mathbb{E}V(X_{n})<\infty$
because the term $b\sum_{k=0}^{n-1}\frac{{\bf 1}_{(X_{k}\in\mathcal{K})}}{\gamma^{k+1}}$
is bounded by $bn$, and so 
\begin{eqnarray*}
E[M_{n+1}] & = & \mathbb{E}\left[\mathbb{E}\left[\frac{V(X_{n+1})}{V_{B}^{*}\gamma^{n+1}}-b\sum_{k=0}^{n}\frac{{\bf 1}_{(X_{k}\in\mathcal{K})}}{\gamma^{k+1}}\mid\mathcal{F}_{n}\right]\right]\\
 & \leq & \mathbb{E}\left[\mathbb{E}\left[\frac{V(X_{n+1})}{V_{B}^{*}\gamma^{n+1}}\mid X_{n}\right]\right]=\mathbb{E}\left[\frac{PV(X_{n})}{V_{B}^{*}\gamma^{n+1}}\right]\\
 & \leq & \mathbb{E}\left[\frac{\gamma V(X_{n})+b}{V_{B}^{*}\gamma^{n+1}}\right]\leq\frac{\gamma\mathbb{E}V(X_{n})+b}{V_{B}^{*}\gamma^{n+1}}<\infty.
\end{eqnarray*}

Therefore, if $\mathbb{E}V(X_{0})<\infty,$ then $\mathbb{E}\mid M_{n}\mid<\infty$
for all $n\geq0,$ and so $M_{n}$ is a super-martingale with the
filtration $\mathcal{F}_{n}$.

Now suppose that $x$ is a fixed point in $B$. Then $1\leq\frac{V(x)}{V_{B}^{*}}$
and ${\bf 1}_{(x\in\mathcal{K})}=0$. Therefore, using the bounded
stopping time $n\wedge\tau_{B}$ for fixed $n$, the optional stopping
theorem implies that 
\begin{align*}
\frac{1}{\gamma^{n}}\mathbb{P}(\tau_{B}>n)  \leq  \mathbb{E}\left[\frac{V(X_{n})}{V_{B}^{*}\gamma^{n}}{\bf 1}_{(\tau_{B}>n)}\right]
 &\leq  \mathbb{E}M_{n\wedge\tau_{B}} \leq  \mathbb{E}M_{0}=\frac{V(x)}{V_{B}^{*}}.
\end{align*}
 Therefore, $\mathbb{P}(\tau_{B}>n)\leq\gamma^{n}\frac{V(x_{0})}{V_{B}^{*}}$,
and so summing over $n$, we have 
\[
\mathbb{E}\tau_{B}\leq\frac{1}{1-\gamma}\cdot\frac{V(x)}{V_{B}^{*}}.
\]
 It follows that if $x\sim\nu$ and $\nu V<\infty$, then $\mathbb{E}\tau_{B}<\infty$. 

Now suppose that $\nu V<\infty$, and define the measure 
\[
\mu_{B}(\cdot)=\mathbb{E}_{\nu}\sum_{k=0}^{\tau_{B}-1}{\bf 1}_{(X_{k}\in\cdot)}.
\]
 We now have that 
\[
\mu_{B}(B)=\mathbb{E}_{\nu}\tau_{B}\leq\frac{1}{1-\gamma}\cdot\frac{\nu V}{V_{B}^{*}}
\]
 so the total mass of $\mu_{B}$ is bounded. Next, let $\xi_{B}(\cdot)=\mathbb{E}{\bf 1}_{(X_{\tau_{B}}\in\cdot)}$
be the measure of the point where $X_{n}$ hits $B$. Using $M_{0}\geq\mathbb{E}M_{n}\geq\mathbb{E}\left[\frac{V(X_{n\wedge\tau_{B}})}{V_{B}^{*}\gamma^{n\wedge\tau_{B}}}\right]{\bf 1}_{\tau_{B}>n}$
and letting $n\rightarrow\infty$, we have 
\[
M_{0}\geq\frac{\mathbb{E}V(X_{\tau_{B}})}{V_{B}^{*}}
\]
 and so
\[
\mathbb{E}V(X_{\tau_{A}})\leq V(X_{0}).
\]
 Therefore, $\xi_{B}V\leq\nu V\leq\infty$. Setting $B=A_{j}\cap\mathcal{K}^{c}$,
and with compactness of $A_{j}\cap\mathcal{K}$, the result follows.
\end{proof}

Our next lemma, which is the key to proving the first main theorem,
relates an injection measure to its occupation and exit measure through
the sub-stochastic kernels. Recall that, when a trajectory to exit
is started, the threshold that determines when an exit occurs, $\kappa$,
is chosen from a distribution $\eta$ on $[0,1]$. 

\begin{lemma}\label{Lemma 2} Suppose that (A0) and (A1) hold. Let
$\nu_{j}$ be a distribution on $A_{j}$ such that $\nu_{j}V<\infty$.
Fix $\kappa\in[0,1]$, and let $\tau_{\kappa,j}$ $\mu_{\kappa,j}$,
$\xi_{\kappa,j}$ be the exit time, occupation, and exit measures
for a trajectory in $A_{j}$ started at $\nu_{j}$, if the threshold
is chosen to be $\kappa$. Then 
\begin{align*}
\mu_{\kappa,j}=\sum_{k=0}^{\infty}\nu_{j}\hat{P}_{\kappa,j}^{k}\ ,\qquad
\mu_{\kappa,j}(I-\hat{P}_{\kappa,j})&=\nu_{j}\ ,\quad\text{and}\quad
\xi_{\kappa,j}=\mu_{\kappa,j}\bar{P}_{j}
\end{align*}
where $I$ is the identity kernel. \end{lemma}

Note that the true occupation and injection measures are obtained
from those above by integrating over $\kappa$: $\mu_{j}(\cdot)=\int_{[0,1]}\mu_{\kappa,j}(\cdot)\eta(d\kappa)$,
$\xi_{j}(\cdot)=\int_{[0,1]}\xi_{\kappa,j}(\cdot)\eta(d\kappa)$. 

\begin{proof} For the first equality, the definition of the occupation
measure gives 
\begin{align*}
 \mu_{\kappa,j}(\cdot)&=\mathbb{E}_{\nu_{j}}\left[\sum_{k=0}^{\tau_{\kappa,{}_{j}}-1}{\bf 1}(X_{k}\in\cdot)\right]\\
\mathbb{E}_{\nu_{j}}\left[\sum_{k=0}^{\infty}{\bf 1}(X_{k}\in\cdot)\cdot{\bf 1}(\tau_{\kappa,\nu_{j}}>k)\right]
&=\sum_{k=0}^{\infty}\mathbb{P}_{\nu_{j}}(X_{k}\in\cdot,\tau_{\kappa,\nu_{j}}>k)=\sum_{k=0}^{\infty}\nu_{j}\hat{P}_{\kappa,j}^{k}(\cdot)
\end{align*}
Because $\nu_{j}\hat{P}_{\kappa,j}^{k}(\cdot)$ is the probability
of a particle being in $\cdot$ at time $k$ and not having left $A_{j}$
at any time before $k.$ The series on the third line above converges
because $\mathbb{E}\tau_{\nu_{j}^{0}}<K$ by assumption (A0). The
second equality is given by applying $I-\hat{P}_{\kappa,j}$ to the
first. 

For the final equality,
\begin{align*}
\xi_{\kappa,j}( \ccdot)=\mathbb{P}_{\nu_{j}}(X_{\tau_{\kappa,j}}\in\ccdot)&=\sum_{k=1}^{\infty}\mathbb{P}_{\nu_{j}}(X_{k}\in\ccdot,\tau_{\kappa,j}=k)=\sum_{k=1}^{\infty}\nu_{j}\hat{P}_{\kappa,j}^{k-1}\bar{P}_{\kappa,j}(\ccdot)\\
&=\left(\sum_{k=0}^{\infty}\nu_{j}\hat{P}_{\kappa,j}^{k}\right)\bar{P}_{\kappa,j}=\mu_{\kappa,j}\bar{P}_{\kappa,j}.
\end{align*}
\end{proof}
We can now put the equalities in the last lemma together to derive,
after integrating over $\kappa$, that a fixed point of $Q$ must
have a corresponding occupation measure which is a fixed point of
$P$, as defined in \ref{Theorem 1}. We then can establish existence
and uniqueness of the fixed point of $Q$, $\nu^{*}$, since $P$
is assumed to have a unique invariant measure. 

\begin{proof}[Proof of Theorem \ref{Theorem 1}] Let $\nu=\sum_{j}a_{j}\nu_{j}$
be a total injection measure, with corresponding occupation measure
$\mu=\sum_{j}a_{j}\mu_{j}.$ By integrating the equations of Lemma~\ref{Lemma 2} 
we get, for each $j$, 
\begin{align*}
\nu_{j}=\int_{[0,1]}\mu_{\kappa,j}(I-\hat{P}_{\kappa,j})\eta(d\kappa)\qquad\text{and}\qquad
\xi_{j}=\int_{[0,1]}\mu_{\kappa,j}\bar{P}_{\kappa,j}\eta(d\kappa).
\end{align*}
 Furthermore, for any $\kappa$ and $j$, $\hat{P}_{\kappa,j}+\bar{P}_{\kappa,j}=P_{\kappa}$
and $\sum_{j}a_{j}\xi_{j}=\nu Q$, by construction. Therefore, we
can subtract the second equation above from the first to get 
\begin{align*}
  \nu-\nu Q=\int_{[0,1]}\sum_{j}a_{j}\mu_{\kappa,j}(I-P_{\kappa})\eta(d\kappa)
=\mu(I-P)
\end{align*}
Therefore, $\nu=\nu Q$, if and only if $\mu(I-P)=0$, in which case
$\mu\sim\pi$, where $\pi$ is the unique invariant distribution of
$P$. 

Now, for a given occupation measure $\mu_{j}$ on $A_{j}$, define
an injection measure $\nu_{j}$ on $A_{j}$ by 
\[
\nu_{j}\sim\mu_{j}(I-\hat{P}_{j})
\]
 recalling that $\hat{P}_{j}=\int_{[0,1]}\hat{P}_{\kappa,j}\eta(d\kappa)$.
Then the occupation measure given by $\nu_{j}$ on $A_{j}$ is, by
Lemma~\ref{Lemma 2}, 
\begin{align*}
\int_{[0,1]}\nu_{j}(I+\hat{P}_{\kappa,j}+\hat{P}_{\kappa,j}^{2}+\cdots)\eta(d\kappa)
&=\int_{[0,1]}\mu_{j}(I-\hat{P}_{j})(I+\hat{P}_{\kappa,j}+\hat{P}_{\kappa,j}^{2}+\cdots)\eta(d\kappa)\\
&=\mu_{j}(I-\hat{P}_{j})(I+\hat{P}_{j}+\hat{P}_{j}^{2}+\cdots)=\mu_{j}
\end{align*} Note that the infinite sum of operators converges, because for any $B\subset A_{j}$ and any probability measure $\eta$, 
\[
\sum_{k=0}^{\infty}\eta \hat{P}_{j}^{k}(B)\leq\sum_{k=0}^{\infty}\eta \hat{P}_{j}^{k}(A_j)=\sum_{k=0}^{\infty}\mathbb{P}_{\eta}(\tau_{j}\geq k)\leq K
\]
Where $K$ is as in (A1). Therefore, given an occupation measure $\mu_{j}$, the associated
injection measure is $\mu_{j}(I-\hat{P}_{j})$, up to normalization. 

We can now show existence of the fixed point. Let $\pi\mid_{A_{j}}$
be the un-normalized restriction of $\pi$ to $A_{j}$. Define the
injection measure on $A_{j}$ by 
\[
\nu_{j}^{*}\sim(\pi\mid_{A_{j}})(I-\hat{P}_{j})
\]
 and the weights by 
\[
a_{j}^{*}=\frac{1}{W}\frac{\pi(A_{j})}{\mathbb{E}\tau_{\nu_{j}^{*}}}
\]
Where $W=\sum_{j}\frac{\pi(A_{j})}{\mathbb{E}\tau_{\nu_{j}^{*}}}.$
Then the above calculation show that if $\nu^{*}=\sum_{j}a_{j}^{*}\nu_{j}^{*}$,
then $\nu^{*}-\nu^{*}Q\sim\pi(I-P)=0$. 

Uniqueness follows similarly. If $\nu Q=\nu$, then we must have $\mu P=\mu$,
so $\mu\sim\pi$. But then $\nu_{j}=\mu_{j}(I-\hat{P}_{j})\sim(\pi\mid_{A_{j}})(I-\hat{P}_{j})$,
and since $\nu_{j}$ must be a probability measure by construction,
this determines all the $\nu_{j}$. Similarly, the weights must be
given by the formula for $a_{j}^{*}$ above, in order to have $\mu\sim\pi$,
so they are determined as well. Therefore, the fixed point of $Q$
is unique. Furthermore, it satisfies $\nu^{*}V<\infty$, because $(\pi\mid_{A_{j}})(I-\hat{P}_{j})V<\infty$
for all $j$. 

\end{proof}

This concludes our proof that the algorithm has the appropriate fixed
distribution. Notice that, in a sense, the argument we used is really
a Poisson equation argument. We use that $\mu_{j}$ satisfies the
Poisson-like equation $\mu_{j}(I-\hat{P}_{j})=\nu_{j}$ for each $j$
to show that, $\mu(I-P)=\nu-\nu Q$, which is like a Poisson equation
on the whole space $A$. This equation is what lets us relate the
fixed point of $P$ to that of $Q$. 

In the rest of the proofs, we approach the problem of showing that
the injection measures given by the algorithm actually approach the
fixed point of $Q$. 

\subsection{Convergence Proof for the Basic Algorithm}

We will first prove that the measures $\nu Q^{n}$, which are the
distributions found by the basic algorithm without finite approximations,
do in fact converge to $\nu^{*}$ as $n\rightarrow\infty$, for any
initial $\nu$. Later we will improve our estimate of how fast the
convergence is. 

The proof we give below has a simple intuition. Consider a particle
following the exit kernel $Q$, whose position at the $n$-th exit
from a stratum is $X_{n}$. Assumption (B1) implies that, after $m$
exits, the particle has at least a certain probability of being in
any stratum. This means that, for two such particles, a coupling can
be constructed so they have a chance, bounded from below, of being
in the same stratum. Assumption (A2) then implies that they have a
probability of being at the same location on the next exit, for an
appropriate coupling. Therefore, the operator that applies the algorithm
$m+1$ times, i.e. $Q^{m+1}$, is a contraction in total variation.
With this intuition, we give the formal proof now. 

\begin{proof}[Proof of Theorem~\ref{Theorem 2}:] First note that,
because $Q$ is a time-homogeneous Markov kernel, it suffices to consider
is a contraction in the case where the initial measures are single-point
delta distributions. Therefore, let $\nu^{1}=x_{1},\nu^{2}=x_{2}$
where $x_{1},x_{2}\in A$ and $x_{1}\neq x_{2},$ and let $n\geq m$. 

Let $j_{1}=\Idx(x_{1})$ and $j_{2}=\Idx(x_{2})$. Let $X_{0}^{i}=x_{i}$,
and $X_{n}^{i}$ be the location of the $n$-th exit from a stratum
starting from $X_{0}^{i},$ for $i=1,2.$ $X_{0}^{i}$ is in $A_{j_{i}},$
so by (A1), with probability at least $c,$ a particle started at
$X_{0}^{i}$ exits as if it started at $\tilde{\nu}_{j_{i}},$ That
is, there exists measures $\eta_{j_{i}}^{i}$ on $A_{j_{i}}^{c}$
for $i=1,2$ such that 
\[
\nu^{i}Q=\delta_{x_{i}}Q=c\tilde{\xi}_{j_{i}}+(1-c)\eta_{j_{i}}^{i}=c\tilde{\nu}_{j_{i}}Q+(1-c)\eta_{j_{i}}^{i}
\]
So applying $Q^{n-1}$ to both sides, where $n\geq m$, 
\[
\nu^{i}Q^{n}=c\tilde{\nu}_{j_{i}}Q^{n}+(1-c)\eta_{j_{i}}^{i}Q^{n-1}.
\]
Now, by assumption (B1), the probability that a particle starting
from the QSD in $A_{j_{i}}$ is in $A_{k}$ after $n$ exits is minorized
by $ua_{k}^{*},$ for any $j$, which means that 
\[
\tilde{\nu}_{j_{i}}Q^{n}(A_{k})\geq ua_{k}^{*}\,\forall i,k.
\]
Therefore, for $i=1,2,$ 
\[
\nu^{i}Q^{n}(A_{k})\geq u\cdot a_{k}^{*}.
\]
Therefore the distributions of $\Idx(X_{n}^{1}),$ $\Idx(X_{n}^{2})$
are both minorized by $cu$ times the same probability vector, $a_{k}^{*}$.
Combining these two observations, there exists a coupling $X_{n}^{1},$ $X_{n}^{2}$ of $\nu^{1}Q^{n}$,
$\nu^{2}Q^{n}$ such that $\mathbb{P}(\Idx(X_{n}^{1})=\Idx(X_{n}^{2}))\geq cu$.
Now, by (A2), conditioned on two particles starting in the same stratum
$A_{k}$, there exists a coupling such that the probability of their
exit points being equal is $\geq c$, because their exit distributions
are both minorized by $c\tilde{\xi}_{k}.$ Therefore, there exists
a coupling $X_{n+1}^{1}$, $X_{n+1}^{2}$ of $\nu^{1}Q^{n+1},$ $\nu^{2}Q^{n+1}$
such that 
\begin{align*}
\mathbb{\mathbb{P}}\big(X_{n+1}^{1}=X_{n+1}^{2}\big) & \geq\mathbb{P}\big(\Idx(X_{n}^{1})=\Idx(X_{n}^{2})\big)\cdot\mathbb{P}\big(X_{n}^{1}=X_{n}^{2}\mid\Idx(X_{n}^{1})=\Idx(X_{n}^{2})\big)\\
 & \geq cu\cdot c=c^{2}u.
\end{align*}

Therefore, $\parallel\nu^{1}Q^{n+1}-\nu^{2}Q^{n+1}\parallel\leq1-c^{2}u=(1-c^{2}u)\parallel\nu^{1}-\nu^{2}\parallel.$

\end{proof}

\subsection{Basic Version - Long Term Convergence Rate}

Proving our more precise estimate of the rate of convergence in the
long term will require more machinery. As we saw above, assumption
(A2) allows us to show that coupling occurs if two particles start
in the same stratum. If the difference between two injection measures,
$\nu^{1}$ and $\nu^{2}$, is mostly in their individual injection
measures (i.e. the $\nu_{j}$), but their weights are mostly the same,
then we can conclude that they get closer when evolved by $Q$. This
is because two particles distributed by $\nu^{1}$, $\nu^{2}$ which
are not already coupled have a good chance of being in the same stratum,
and then (A2) says they have a chance of coupling on the next exit.
But what if the difference is mostly in the weights? Then the particles
have little chance to couple on the next exit, since they are most
likely not in the same stratum. So we can't conclude that $Q$ is
a contraction. This can be thought of as the ``problem case'', as
far as a coupling argument is concerned. 

We solve this issue with two tools. First, given the equilibrium $\nu^{*}$
and a given measure $\nu$, we want to decompose them into the the
two above cases. We can start with the standard way of decomposing
the measures into equal and mutually singular parts: 
\[
\nu=(1-\epsilon)\bar{\nu}+\epsilon\nu^{1}
\]
\[
\nu^{*}=(1-\epsilon)\bar{\nu}+\epsilon\nu^{2}
\]
 where $\bar{\nu}$, $\nu^1$, and $\nu^2$ are probability measures with $\nu^{1}\perp\nu^{2}$ and $\epsilon=\parallel\nu-\nu^{*}\parallel$.
We can further decompose the mutually singular parts into a part that
looks like the ``easy case'' above, and a part that looks like the
``problem case''. We get the following expression: 
\begin{equation}\label{eq:2}
  \begin{aligned}
    \nu^{1} & =(1-\gamma)\sum_{j}b_{j}^{1}\nu_{j}+\gamma\sum_{j}p_{j}\eta_{j}^{1}\\
\nu^{2} & =(1-\gamma)\sum_{j}b_{j}^{2}\nu_{j}^{*}+\gamma\sum_{j}p_{j}\eta_{j}^{2}
  \end{aligned}
\end{equation}
 for some probability vectors $p,$ $b^{1},$ $b^{2}$ and measures
$\eta^{1},\eta^{2}$ where $b^{1}\perp b^{2}$, $\eta^{1}\perp\eta^{2}$,
and some $\gamma\in[0,1]$. The first term represents the parts of
$\nu^{1},\nu^{2}$ that have different weights. The second term represents
the parts that have the same weights but orthogonal distributions
on the strata. $\gamma$ is the share of the second term relative
to all of the TV distance between $\nu$ and $\nu^{*}$. We will handle
the cases where $\gamma$ is large or small separately, and get a
bound for each. 

However, as explained above, $Q$ may not actually bring the two measures
closer together in the ``problem case'', where $\gamma<<1$. To
fix this, we introduce the following new metric, which is equivalent
to TV distance in the sense that they induce the same topology: 
\begin{equation*}
  d_{\alpha}(\nu,\nu^{*})=\parallel\nu-\nu^{*}\parallel_{TV(A)}+\alpha d_{w}(a,a^{*})
\end{equation*}
 where $a,a^{*}$ are the weight vectors of $\nu,\nu^{*}$, and $G^{*}$
is the transition matrix between strata in equilibrium. The metric
$d_{w}$ on weight vectors is given by 
\[
d_{w}(a,b)=\sum_{k=0}^{m}\lambda^{k}\parallel a(G^{*})^{k}-b(G^{*})^{k}\parallel+\sum_{k=m+1}^{\infty}\lambda^{m}\parallel a(G^{*})^{k}-b(G^{*})^{k}\parallel
\]
 where $\lambda,m$ are as in (B2). The intuition behind the new metric
is as follows. Even if we are in the problem case, the weights of
our injection measure should get closer to the true weights, which
should allow for coupling to occur at some later step. So we add a
term to the standard TV-distance that should contract if the weights
move approximately by $G^{*}$ under $Q$. We choose the metric $d_{w}$
on probability vectors because under it, $G^{*}$ is a contraction
with constant $\lambda$. 

With the decomposition and metric defined, our strategy is now as
follows: First, we consider the ``easy case'', where $\gamma$ is
close to $1$. Then (A1) lets us conclude that coupling occurs on
the next exit, so $\nu$ gets closer to $\nu^{*}$ when acted on by
$Q$. In the ``problem case'', where $\gamma$ is close to $0$,
the weights of $\nu$ get closer to those of $\nu^{*}$ because they
move approximately by $G^{*},$ so we still get contraction in $d_{\alpha}$.
The complication we encounter is that the weights do not move exactly
by $G^{*}$, but by the transition matrix $G$ given by $\nu$. The
result is that we only get local contraction of $Q$ in $d_{\alpha}$.
But we already have global convergence by \ref{Theorem 2}. So we
still get a rate of convergence that eventually applies. 

With the strategy now laid out, we can now proceed through the proof
of \ref{Theorem 3}. We begin by showing the decomposition (\ref{eq:2}). 

\begin{lemma}\label{Lemma 3}Let $\nu=\sum a_{j}\nu_{j}$ be an injection
measure on $A$, and let $\nu^{*}=\sum_{j}a_{j}^{*}\nu_{j}^{*}$ be
the fixed point of $Q$ then there exist probability vectors $p,$
$b^{1}\perp b^{2}$, measures $\eta^{1}\perp\eta^{2}$ and $\gamma\in[0,1]$
such that (\ref{eq:2}) holds. 

\end{lemma}

\begin{proof}Start by decomposing the weight vectors of $\nu$, $\nu^{*}$
into equal and orthogonal parts: 
\[
a=(1-\hat{\epsilon})\bar{a}+\hat{\epsilon}b^{1}
\]
\[
a^{*}=(1-\hat{\epsilon})\bar{a}+\hat{\epsilon}b^{2}.
\]
with $b^{1}\perp b^{2}$. Similarly, we can decompose the injection
measures on each stratum: 
\[
\nu_{j}=(1-\epsilon_{j})\bar{\nu}_{j}+\epsilon_{j}\eta_{j}^{1}
\]
\[
\nu_{j}^{*}=(1-\epsilon_{j})\bar{\nu}_{j}+\epsilon_{j}\eta_{j}^{2}.
\]
 We can now substitute these expressions for the injection measures
and weights into $\nu=\sum a_{j}\nu_{j}$, $\nu^{*}=\sum_{j}a_{j}^{*}\nu_{j}^{*}$.
We get 
\begin{align*}
\nu^{0} & =(1-\hat{\epsilon})(1-Z)\left[\sum_{j}\frac{\bar{a}_{j}(1-\epsilon_{j})}{1-Z}\bar{\nu}_{j}\right]+\hat{\epsilon}\left[\sum_{j}b_{j}^{1}\nu_{j}\right]+(1-\hat{\epsilon})Z\left[\sum_{j}p_{j}\eta_{j}^{1}\right]\\
\nu^{*} & =(1-\hat{\epsilon})(1-Z)\left[\sum_{j}\frac{\bar{a}_{j}(1-\epsilon_{j})}{1-Z}\bar{\nu}_{j}\right]+\hat{\epsilon}\left[\sum_{j}b_{j}^{*}\nu_{j}^{*}\right]+(1-\hat{\epsilon})Z\left[\sum_{j}p_{j}\eta_{j}^{2}\right]
\end{align*}
where $Z=\sum_{j}\bar{a}_{j}\epsilon_{j}$ and $p_{j}=\frac{\bar{a}_{j}\epsilon_{j}}{Z}$.
Setting $\gamma=\frac{(1-\hat{\epsilon})Z}{\hat{\epsilon}+(1-\hat{\epsilon})Z},$the
result follows. 

\end{proof}

As mentioned above, our proof strategy is complicated by the fact
that the weight vectors of two injection measures move by different
matrices under $Q$. However, if all $a_{j}^{*}$ are non-zero, i.e.
every stratum has some weight in equilibrium, then $G,G^{*}$ will
be close if $\nu,\nu^{*}$ are close in total variation. So our next
step is to bound the difference in weight vectors after applying $Q,$
in the case where the transition matrices are close. 

\begin{lemma} \label{Lemma 4}Suppose that (A0)-(A2) and (B2) hold.
If $G$, $G^{*}$ are the transition matrices given by the exit measures
for $\nu$, $\nu^{*}$, and if $\parallel G-G^{*}\parallel_{\infty}<\delta$,
then 
\begin{equation}
d_{w}(aG,a^{*}G^{*})\leq S\left((1-c)\gamma+\delta(1-\gamma)\right)\parallel\nu-\nu^{*}\parallel_{TV}+\lambda d_{w}(a,a^{*}).\label{eq:3}
\end{equation}
 Where $S=1+\lambda+\lambda^{2}+\lambda^{3}+\cdots=\frac{1}{1-\lambda}$.
\end{lemma}

\begin{proof}

By decomposing the weight vectors $a$, $a^{*}$ into equal and orthogonal
parts, as in the proof of Lemma~\ref{Lemma 3}, we have 
\begin{eqnarray*}
aG-a^{*}G^{*} & = & a(G-G^{*})+(a-a^{*})G^{*}\\
 & = & (1-\hat{\epsilon})\bar{a}(G-G^{*})+\hat{\epsilon}b^{1}(G-G^{*})+(a-a^{*})G^{*}.
\end{eqnarray*}
 Therefore, by the triangle inequality, 
\begin{equation*}
d_{w}(aG,a^{*}G^{*}) \leq  (1-\hat{\epsilon})d_{w}(\bar{a}G,\bar{a}G^{*})+  \hat{\epsilon}d_{w}(aG^{*},a^{*}G^{*})  +  \hat{\epsilon}d_{w}(b^{1}G,b^{1}G^{*}).
\end{equation*}
 Now, $G^{*}$ is a contraction in $d_{w}$, so $d_{w}(a^{0}G^{*},a^{*}G^{*})\leq\lambda d_{w}(a^{0},a^{*})$.
The construction of $d_{w}$ and assumption (B2) imply that 
\[
d_{w}(b^{1}G,bG^{*})\leq S\parallel b^{1}G-b^{1}G^{*}\parallel
\]
which is $\leq S\parallel G-G^{*}\parallel_{\infty}\leq s\delta$
because $b^{1}$ is a probability vector. 

Using the definition of $G,$ $G^{*}$ and assumption (A2), 
\begin{align*}
\Big|\left(\bar{a}(G-G^{*})\right)_{j}\Big| & =  \Big|\sum_{k}\bar{a}_{k}(G_{jk}-G_{jk}^{*}\Big| =  \Big|\sum_{k}\bar{a}_{k}(\xi_{k}(A_{j})-\xi_{k}^{*}(A_{j})) \Big|\\
 & \leq \sum_{k}\bar{a}_{k}\parallel\xi_{j}-\xi_{j}^{*}\parallel \leq  (1-c)\sum_{k}\bar{a}_{k}\parallel\nu_{j}-\nu_{j}^{*}\parallel_{TV}.
\end{align*}
Therefore, 
\begin{eqnarray*}
(1-\hat{\epsilon})\parallel\bar{a}(G-G^{*})(G^{*})^{k}\parallel_{TV} & \leq & (1-c)(1-\hat{\epsilon})\sum\bar{a}_{k}\epsilon_{k}\\
 & = & (1-c)\gamma\parallel\nu-\nu^{*}\parallel_{TV}
\end{eqnarray*}
 Which implies that $(1-\hat{\epsilon})d_{w}(\bar{a}G,\bar{a}G^{*})\leq(1-c)\gamma S\parallel\nu-\nu^{*}\parallel_{TV}$.
Putting the inequalities just derived together gives the desired result.
\end{proof}

Now we are ready to show that $Q$ acts as a local contraction in
$d_{\alpha}$, for both small and large $\gamma$. We start with the
case where $\gamma$ is bounded away from $0$. In this case, the
idea is that, if two particles coupling $\nu$ and $\nu^{*}$ are
not at equal positions, they have a chance, bounded below, of being
in the same stratum. (A2) then lets us assure that they can couple
when moving by $Q$. The rest of the proof of the following lemma
is simply keeping track of how all the terms in $d_{\alpha}(\nu Q,\nu^{*}Q)$
behave. 

\begin{lemma} \label{Lemma 5}Suppose that (A0)-(A2) and (B2) hold
and fix $\beta\in(0,1)$. Then for initial $\nu$ and equilibrium
$\nu$ injection measures, if $\gamma>\beta$ and $\parallel G-G^{*}\parallel_{\infty}<\delta$,
\[
d_{\alpha}(\nu Q,\nu^{*}Q)\leq\max\left(1-\beta c+\alpha S(1-c)+\alpha S\delta,\lambda\right)\cdot d_{\alpha}(\nu,\nu^{*}).
\]
\end{lemma}

\begin{proof}

Suppose $\nu\neq\nu^{*}$, since otherwise the result is trivial.
Let $X,X^{*}$ be a coupling of $\nu,\nu^{*}$. By Lemma~\ref{Lemma 3}, the coupling
can be chosen so that with probability $\gamma\parallel\nu-\nu^{*}\parallel$,
$\Idx(X)=\Idx(X^{*})$ and $X\neq X^{*}$, so that 
\begin{eqnarray*}
\mathbb{P}\Big(\Idx(X)=\Idx(X^{*})\mid X\neq X^{*}\Big) & = & \frac{\mathbb{P}\Big(\Idx(X)=\Idx(X^{*}),X\neq X^{*}\Big)}{\mathbb{P}(X=X^{*})}\\
 & = & \frac{\gamma\parallel\nu-\nu^{*}\parallel}{\parallel\nu-\nu^{*}\parallel}=\gamma\\
 & \geq & \beta.
\end{eqnarray*}
 Assumption (A2) now implies that if $Y,Y^{*}$ are a coupling of
$\nu Q,\nu^{*}Q=\nu^{*}$, then they can be chosen so that 
\[
\mathbb{P}(Y=Y^{*}\mid X\neq X^{*})\geq c\gamma\geq c\beta.
\]
 Therefore, 
\[
\parallel\nu Q-\nu^{*}\parallel\leq(1-c\beta)\parallel\nu-\nu^{*}\parallel.
\]
 Combining this with Lemma~\ref{Lemma 4} and the definition of $d_{\alpha}$
implies that 
\[
d_{\alpha}(\nu Q,\nu Q)\leq\left(1-\beta c+\alpha S(1-c)+\alpha S\delta\right)\parallel\nu-\nu^{*}\parallel_{TV}+\alpha\lambda d_{w}(a,a^{*})
\]
\[
\leq\max\left(1-\beta c+\alpha S(1-c)+\alpha S\parallel G^{0}-G^{*}\parallel_{\infty},\lambda\right)d_{\alpha}(\nu^{0},\nu^{*}).
\]

\end{proof}

Now we can approach the case where $\gamma$ is bounded away from
$1$ in a similar way. The idea this time is that, if we are close
to the problem case, but $\parallel G-G^{*}\parallel$ is small, then
the weights will get more accurate under $Q$ and so the distance
in $d_{\alpha}$ should still contract, even if no coupling occurs. 

\begin{lemma} \label{Lemma 6} Suppose that (B2) holds, and fix $\beta\in(0,1)$.
Suppose that $\alpha<\frac{1}{S(1-c)}$. Then for any $\nu$, $\nu^{*}$
such that $\gamma\leq\beta$, $\parallel G-G^{*}\parallel<\delta$,
\begin{equation}
d_{\alpha}(\nu Q,\nu^{*}Q)\leq\frac{\left(1+\alpha S\delta\right)+(1-\beta)\alpha\lambda}{1+(1-\beta)\alpha}d_{\alpha}(\nu,\nu^{*}).\label{eq:4}
\end{equation}
\end{lemma}

\begin{proof} 

As before, we can assume that $\nu\neq\nu^{*}$. First, as in the
previous proof, we have 
\[
\parallel\nu Q-\nu^{*}Q\parallel\leq(1-\gamma c)\parallel\nu-\nu^{*}\parallel
\]
 and by the construction of $\gamma$, 
\begin{align*}
\parallel a-a^{*}\parallel & =  (1-\gamma)\parallel\nu-\nu^{*}\parallel \geq  (1-\beta)\parallel\nu-\nu^{*}\parallel\\
\Longrightarrow&\parallel\nu-\nu^{*}\parallel  \leq  \frac{1}{1-\beta}\parallel a^{0}-a^{*}\parallel\leq \frac{1}{1-\beta}d(a,a^{*}).
\end{align*}

Now, putting the above inequalities together with the result of Lemma~\ref{Lemma 4}
implies that 
\begin{align*}
d_{\alpha}(\nu,\nu^{*}) &\leq  \big(1+(\alpha S(1-c)-c)\gamma+\alpha S\delta\big)\parallel\nu^{0}-\nu^{*}\parallel_{TV}+\lambda d_{w}(a^{0},a^{*})\\
 & \leq \rho^{0}d_{\alpha}(\nu^{0},\nu^{*})
\end{align*}
 where 
\[
\rho=\frac{\left(1+(\alpha S(1-c)-c)\gamma_{0}+\alpha S\delta\right)\parallel\nu-\nu^{*}\parallel+\alpha\lambda d_{w}(a,a^{*})}{\parallel\nu-\nu^{*}\parallel+d_{w}(a,a^{*})}.
\]
 Simplifying $\rho$ gives the desired result. 

\end{proof}

Notice that the contraction constants found in the previous two lemmas
are $<1$ if $\delta$ is sufficiently small. This is why we can conclude
that $Q$ is a local contraction in $d_{\alpha}$. With both cases
now out of the way, we can put them together to get our next main
convergence theorem. 

\begin{proof}[Proof of Theorem~\ref{Theorem 3}:] Theorem~\ref{Theorem 2}
implies that, under the given assumptions, $\parallel\nu^{n}-\nu^{*}\parallel\rightarrow0$,
where $\nu^{n}=\nu Q^{n}$. Since, by assumption, $a_{j}^{*}>0$ for
all $j$, this implies that $\nu_{j}^{n}\rightarrow\nu_{j}^{*}$ for
all $j,$ in total variation, and therefore $G^{n}\rightarrow G^{*}$,
where $G^{n}$ is the transition matrix given by $\nu^{n}$. Therefore,
there exists $l\geq0,\delta>0$ such that, for all $n\geq l$ and
$\alpha<\frac{c}{S(1-c)}$, $\parallel G-G^{*}\parallel_{\infty}<\delta$
and 
\[
q_{\alpha,\beta}^{\delta}:=\max\left(1-\beta c+\alpha S(1-c)+\alpha S\delta,\lambda,\frac{\left(1+\alpha S\delta\right)+(1-\beta)\alpha\lambda}{1+(1-\beta)\alpha}\right)<1.
\]
 It follows that $d_{\alpha}(\nu^{n+1},\nu^{*})\leq q_{\alpha,\beta}^{\delta}\cdot d_{\alpha}(\nu^{n},\nu^{*})$
for $n\geq l.$ Furthermore, as $l\rightarrow\infty$, the above holds
for arbitrarily small $\delta$. Also note that for $\lambda<1$,
$\frac{\left(1+\alpha S\delta\right)+(1-\beta)\alpha\lambda}{1+(1-\beta)\alpha}\geq\frac{1+(1-\beta)\alpha\lambda}{1+(1-\beta)\alpha}>\lambda$.
Therefore, there exist $q_{n}\in(0,1)$ such that $d_{\alpha}(\nu^{n+1},\nu^{*})\leq q^{n}\cdot d_{\alpha}(\nu^{n},\nu^{*})$
for sufficiently large $n$, and 
\[
q^{n}\rightarrow\lim_{\delta\rightarrow0}q_{\alpha,\beta}^{\delta}=\max\left(1-\beta c+\alpha S(1-c),\frac{1+(1-\beta)\alpha\lambda}{1+(1-\beta)\alpha}\right).
\]
 Finally, optimizing the limit of $q^{n}$ over all choices of $\alpha$,
$\beta$ proves the desired result. 

\end{proof}

\subsection{Local Convergence of the Eigenvector Version}

Now we turn to proving our final result, that the eigenvector version
converges locally with a rate we can bound. The strategy again revolves
around the decomposition (\ref{eq:2}). This time, however, we do
not have to deal with the ``problem case''. The reason is that,
as a result of the perturbation bound in (B3), we can show that the
parameter $\gamma$ is never too small when the weights are given
by the eigenvector of the transition matrix $G$ associated to $\{\nu_{j}\}$.
Therefore, we always get that coupling can occur, by assumption (A1). 

However, we also run into a new problem: this version of the algorithm
cannot be represented as evolving the measure $\nu$ forward by a
single Markov kernel, the way the basic version evolves by $Q$. The
reason is that each of the eigenvector weights depends on all the
measures $\nu_{j}$. The eigenvector version does not act meaningfully
on just one $\nu_{j}$. Therefore, we cannot assure that no decoupling
occurs in the step where the new weights are chosen. However, similarly
to the previous proof, we can control the size of this decoupling
in the case where the $\nu_{j}$ are already close to $\nu_{j}^{*}$.
Therefore, we will end up with a local convergence result. 

Our first goal is to show that, when using the eigenvector weights,
the parameter $\gamma$ is bounded below. This follows somewhat straightforwardly
from assumption (B3), and the fact that the transition matrix entries
are themselves determined by the injection measures. Note that the
eigenvector weights for the equilibrium injection measure are simply
given by $a^{*}$, because by construction of the fixed point of $Q,$
$a^{*}G^{*}=a^{*}$. 

\begin{lemma} \label{Lemma 7} Suppose that (B3) holds. Let $G^{*}$
be the transition matrix given by the fixed point measures $\nu_{j}$,
with principal eigenvector $a^{*}$, i.e. $a^{*}G^{*}=a^{*}$. Let
$G$ be another transition matrix on $\{1,\ldots,J\}$ with eigenvector
$z$, such that $z_{j}>0$ for all $j$. Then for any $\delta>0$,
there exist constants $\theta_{ik}$ for $i\neq k$ such that 
\[
\frac{\mid z_{j}-a_{j}^{*}\mid}{z_{j}\wedge a_{j}^{*}}\leq\sum_{i\neq k}\theta_{ik}^{\prime}\mid G_{ik}-G_{ik}^{*}\mid.
\]
if $\sup_{i,k}\mid G_{ik}-G_{ik}^{*}\mid<\delta$. Furthermore, as
$\delta\rightarrow0$, $\theta_{ik}^{\prime}\rightarrow e^{\theta_{ik}}-1$.
\end{lemma}

\begin{proof}

By (B3), there exist $\theta_{ik}>0$ $\forall i,k\in\{1,...,J\}$
such that for all $j\leq J,$ 
\begin{multline*}
  \log\left(\max\left(\frac{z_{j}}{a_{j}^{*}},\frac{a_{j}^{*}}{z_{j}}\right)\right)\leq\sum_{i\neq k}\theta_{ik}\mid G_{ik}-G_{ik}^{*}\mid\\
\Longrightarrow\frac{z_{j}\vee a_{j}^{*}}{z_{j}\wedge a_{j}^{*}}=1+\frac{\mid z_{j}-a_{j}^{*}\mid}{z_{j}\wedge a_{j}^{*}}\leq\exp\Big(\sum_{i\neq k}\theta_{ik}\mid G_{ik}-G_{ik}^{*}\mid\Big).
\end{multline*}
Now, for any $M>0,$ there exists $K>0$ such that $e^{x}-1\leq Kx$
for $0\leq x\leq M.$ Specifically, since $\mid G_{ik}-G_{ik}^{*}\mid\leq1,$
we have 
\begin{align*}
  \exp\Big(\sum_{i\neq k}\theta_{ik}\mid G_{ik}-G_{ik}^{*}\mid\Big)&-=\prod_{i\neq k}\Big(1+\sum_{m=1}^{\infty}\frac{(\theta_{ik}\mid G_{ik}-G_{ik}^{*}\mid)^{m}}{m!}\Big)-1\\
&\leq\prod_{i\neq k}\Big(1+\mid G_{ik}-G_{ik}^{*}\mid\sum_{m=1}^{\infty}\frac{\theta_{ik}^{m}}{m!}\Big)-1=\prod_{i\neq k}\Big(1+\mid G_{ik}-G_{ik}^{*}\mid(e^{\theta_{ik}}-1)\Big)-1\\
&\leq\sum_{i\neq k}\mid G_{ik}-G_{ik}^{*}\mid\left((e^{\theta_{ik}}-1)+O(\mid G-G^{*}\mid)\right),
\end{align*}
where $O(\mid G-G^{*}\mid)$ represents terms which are order $1$
or higher in $\{\mid G_{i^{\prime}k^{\prime}}-G_{i^{\prime}k^{\prime}}^{*}\mid\}_{i\neq k}.$
Since there are only finitely many of these terms, setting $\theta_{ik}^{\prime}$
to be the coefficient of $\mid G_{ik}-G_{ik}^{*}\mid$ in
\begin{equation*}
  \prod_{i\neq k}\left(1+\mid G_{ik}-G_{ik}^{*}\mid(e^{\theta_{ik}}-1)\right)-1\,,
\end{equation*}
with all terms of the form $\mid G_{i^{\prime}k^{\prime}}-G_{i^{\prime}k^{\prime}}^{*}\mid$
replaced by $1$, suffices.

\end{proof}

Now we can use the above lemma to show that $\gamma$ is bounded in
the eigenvector version. 

\begin{lemma} \label{Lemma 8} Suppose that (A0)-(A2) and (B3) hold.
Let $\nu_{j}$ be injection measures for $1,\ldots,J$, with transition
matrix $G$ and corresponding (normalized) eigenvector $z$. Let $\nu=\sum_{j}z_{j}\nu_{j}$,
and decompose $\nu$, $\nu^{*}$ as in (\ref{eq:2}). Then there exists
$r$, independent of $\nu_{j}$, such that 
\begin{equation}
\gamma\geq\frac{1}{1+r}.\label{eq:5}
\end{equation}
 Furthermore, as $\nu_{j}\rightarrow\nu_{j}^{*}$ in TV, $r$ can
be chosen to be arbitrarily close to $r^{\infty}$, as in Theorem~\ref{Theorem 4}. 

\end{lemma}

\begin{proof} 

First, note that by (A2), $G^{*},G\geq c\tilde{G}$, where $\tilde{G}$
is the transition matrix if the injection measure in $A_{j}$ is the
$QSD$ $\tilde{\nu}_{j}$. Therefore, (B3) applies to $G$ and $G^{*}$.
Now write $\hat{\epsilon}=\parallel z-a^{*}\parallel$ and $\epsilon_{j}=\parallel\nu_{j}-\nu_{j}^{*}\parallel$
for all $j.$ 

From Lemma~\ref{Lemma 7}, we have that for all $j,$ 
\begin{align*}
  \frac{\mid z_{j}-a_{j}^{*}\mid}{z_{j}\wedge
  a_{j}^{*}}\leq\sum_{i\neq k}\theta_{ik}^{\prime}\mid
  G_{ik}-G_{ik}^{*}\mid&=\sum_{i}\left(\sum_{k\neq
  i}\theta_{ik}^{\prime}\mid G_{ik}-G_{ik}^{*}\mid\right)\\
  &\leq\sum_{i}\left(\sup_{k}(\theta_{ik}^{\prime})\sum_{k\neq i}\mid
    G_{ik}-G_{ik}^{*}\mid\right)\leq\sum_{i}\sup_{k}(\theta_{ik}^{\prime})\cdot2\parallel{\bf
    G}_{i}-{\bf G}_{i}^{*}\parallel\\
  &\leq\sum_{i}2\cdot\sup_{k}(\theta_{ik}^{\prime})\cdot(1-c)\epsilon_{i}\leq2(1-c)\sup_{i,k}\left(\frac{\theta_{ik}^{\prime}}{a_{i}^{*}}\right)\sum_{i}a_{i}^{*}\epsilon_{i}
\end{align*}
Where the third line follows from (A2), which implies that $\parallel{\bf G}_{i}-{\bf G}_{i}^{*}\parallel\leq\parallel\nu_{j}Q-\nu_{j}^{*}Q\parallel\leq(1-c)\epsilon_{j}$.
Now, (B3) implies that 
\begin{align*}
  \frac{a_{i}^{*}}{z_{i}\wedge a_{i}^{*}}\leq\frac{z_{i}\vee a_{i}^{*}}{z_{i}\wedge a_{i}^{*}}&\leq\exp\left(\sum_{i\neq k}\theta_{ik}\mid G_{ik}-G_{ik}^{*}\mid\right)\\
&\leq\exp\left(2(1-c)\sum_{i}\sup_{k\neq i}(\theta_{ik})\right)=:E
\end{align*}
Therefore, for all $j\in\{1,...,J\},$ 
\begin{align*}
  \hat{\epsilon}\leq\mid z_{j}-a_{j}^{*}\mid&\leq(z_{j}\wedge a_{j}^{*})\cdot2(1-c)\sup_{i,k}\left(\frac{\theta_{ik}^{\prime}}{a_{i}^{*}}\right)E\sum_{i}(z_{i}\wedge a_{i}^{*})\epsilon_{i}\\
&\leq r\sum_{i}(z_{i}\wedge a_{i}^{*})\epsilon_{i}
\end{align*}
where
\[
r=2(1-c)\sup_{j}(a_{j}^{*})\sup_{i,k}\left(\frac{\theta_{ik}^{\prime}}{a_{i}^{*}}\right)E.
\]
 We also now have that $r\rightarrow r^{\infty}$ as $\nu\rightarrow\nu^{*}$,
because $\theta_{ik}^{\prime}\rightarrow e^{\theta_{ik}}-1$. Now, observe
that in the notation of Lemma~\ref{Lemma 3}, 
\[
\frac{(1-\hat{\epsilon})Z}{(1-\hat{\epsilon})Z+\hat{\epsilon}}=\frac{\sum_{j}(z_{j}\wedge a_{j}^{*})\epsilon_{j}}{\sum_{j}(z_{j}\wedge a_{j}^{*})\epsilon_{j}+\hat{\epsilon}}
\]
 and so the result follows. 

\end{proof}

We can now show that the eigenvector version acts as a local contraction
on injection measures. We choose to state the local contraction in
another metric, for which the calculations are simpler. Given injection
measures $\nu_{j}$, with eigenvector weights $z_{j}$ and $\nu=\sum_{j}z_{j}\nu_{j}$,
define 
\[
d(\nu,\nu^{*})=\sum_{j}a_{j}^{*}\parallel\nu_{j}-\nu_{j}^{*}\parallel=\sum_{j}a_{j}^{*}\epsilon_{j}.
\]
 Note that, since the measures $\nu_{j}$ determine the weights $z_{j}$,
this is a well-defined metric on the set of injection measures with
eigenvector weights. The proof of Lemma~\ref{Lemma 8} also implies that $\hat{\epsilon}<rd(\nu,\nu^{*})$
so $d$ is equivalent to total variation distance for such choices
of $\nu$. Since we are in the case where the measures determine the
weights, $d$ can be thought of as the metric that only looks at the
difference in measures, appropriately weighted. By moving between
$d(\nu,\nu^{*})$ and $\parallel\nu-\nu^{*}\parallel$, we can use
a coupling argument and the boundedness of $\gamma$ to show local
contraction in $d$, as we do below. 

\begin{lemma}\label{Lemma 9} Suppose that (A0)-(A2)$.$ Let $\nu^{\prime}=\sum_{j}z_{j}^{\prime}\nu_{j}^{\prime}$
be the measure obtained from $\nu$ via the eigenvector injection
measure method, with $\epsilon_{j}^{\prime}=\parallel\nu_{j}^{\prime}-\nu_{j}^{*}\parallel$.
Then
\[
\sum_{j}(z_{j}\wedge a_{j}^{*})\epsilon_{j}^{\prime}\leq q\sum_{j}(z_{j}\wedge a_{j}^{*})\epsilon_{j}
\]

Where 
\[
q:=1-\frac{1}{1+r}c
\]
\end{lemma}

\begin{proof} First, by Lemma~\ref{Lemma 8} and the same coupling argument as
in Lemma~\ref{Lemma 5}, 

\begin{equation}
\parallel\nu Q-\nu^{*}\parallel\leq(1-\gamma c)\parallel\nu-\nu^{*}\parallel\leq q\parallel\nu-\nu^{*}\parallel.\label{eq:6}
\end{equation}
 The individual injection measures of $\nu Q$ are $\nu_{j}^{\prime}$,
by construction of the algorithm. Furthermore, the weights of $\nu Q$
are $z,$ because $zG=z.$ Therefore, 
\[
\parallel\nu Q-\nu^{*}\parallel=\sum_{j}(z_{j}\wedge a_{j}^{*})\epsilon_{j}+\hat{\epsilon}.
\]
 We also have 
\[
\parallel\nu-\nu^{*}\parallel=\sum_{j}(z_{j}\wedge a_{j}^{*})\epsilon_{j}^{\prime}+\hat{\epsilon}.
\]
 Substituting the above expressions into (\ref{eq:6}), we get 
\[
\sum_{j}(z_{j}\wedge a_{j}^{*})\epsilon_{j}^{\prime}\leq q\sum_{j}(z_{j}\wedge a_{j}^{*})\epsilon_{j}-(1-q)\hat{\epsilon}\leq q\sum_{j}(z_{j}\wedge a_{j}^{*})\epsilon_{j}.
\]

\end{proof} 

$\,$

\begin{lemma} \label{Lemma 10} Fix $\delta>0$. Under the assumptions
of Lemma~\ref{Lemma 9}, if $2(1-c)(\sup_{i\neq k}\frac{\theta_{ik}}{a_{i}^{*}})d(\nu,\nu^{*})<\delta$,
then
\[
d(\nu^{\prime},\nu^{*})\leq q(1+\delta)d(\nu,\nu^{*})
\]

\end{lemma}

\begin{proof} Continuing on from Lemma~\ref{Lemma 9}, we have 
  \begin{align*}
   d(\nu^{\prime},\nu^{*})&=\sum_{j}a_{j}^{*}\epsilon_{j}^{\prime}\leq\sup_{j}\left(\frac{z_{j}\vee a_{j}^{*}}{z_{j}\wedge a_{j}^{*}}\right)\sum_{j}(z_{j}\wedge a_{j}^{*})\epsilon_{j}^{\prime}
\leq q\frac{z_{j}\vee a_{j}^{*}}{z_{j}\wedge a_{j}^{*}}\left(\sum_{j}(z_{j}\wedge a_{j}^{*})\epsilon_{j}\right)\\
&\leq q\cdot\left(1+\frac{\mid z_{j}-a_{j}^{*}\mid}{z_{j}\wedge a_{j}^{*}}\right)\sum_{j}a_{j}^{*}\epsilon_{j}\\
&\leq q\cdot\exp\left(1+2(1-c)\sum_{i}(\sup_{k\neq i}\theta_{ik})\epsilon_{i}\right)d(\nu,\nu^{*})\leq q(1+\delta)d(\nu,\nu^{*}) 
  \end{align*}
 where the second line holds by (A2), and the last line holds because
$d(\nu,\nu^{*})\geq\frac{1}{\sup_{i\neq k}\frac{\theta_{ik}}{a_{i}^{*}}}\cdot\sum_{i}\sup_{k\neq i}(\theta_{ik})\epsilon_{i}$. 

\end{proof}

Note, in particular, that if $\delta<\frac{1}{q}-1$, then $q(1+\delta)<1$,
and so have have local contraction. It is now straightforward to use
the equivalence of TV distance and $d$, when the weights are given
by the eigenvector of $G$, to prove our final convergence result.

\begin{proof}[Proof of Theorem~\ref{Theorem 4}:]

Suppose that the initial injection measure is $\nu^{0}=\sum_{j}z_{j}^{0}\nu_{j}^{0}$,
and that successive $\nu^{n}$ are given by steps of the eigenvector
version, with. Lemma~\ref{Lemma 10} implies that, if $2(1-c)(\sup_{i\neq k}\frac{\theta_{ik}}{a_{i}^{*}})d(\nu^{0},\nu^{*})<\frac{1}{q}-1$,
then $d(\nu^{n},\nu^{*})\rightarrow0$, and for all $n$, $d(\nu^{n+1},\nu^{*})\leq p^{n}d(\nu^{n},\nu^{*})$
where $p^{n}$ is the value $q(1+\delta)$, computed from $\nu=\nu^{n}$.
By construction of $p^{n}$, we have 
\[
p^{n}\rightarrow1-\frac{c}{1+r^{\infty}}.
\]
 Next, by the bounds derived in Lemma~\ref{Lemma 8}, 
 \begin{align*}
   \parallel\nu^{0}-\nu^{*}\parallel=\parallel z^{0}-z^{*}\parallel+\sum_{j}(z_{j}^{0}\wedge z_{j}^{*})\parallel\nu_{j}^{0}-\nu_{j}^{*}\parallel&\leq(1+r)\sum_{j}(z_{j}^{0}\wedge z_{j}^{*})\parallel\nu_{j}^{0}-\nu_{j}^{*}\parallel\\
&\leq(1+r)d(\nu^{0},\nu^{*}).
 \end{align*}
 Therefore, we have convergence if $\frac{2(1-c)(\sup_{i\neq k}\frac{\theta_{ik}}{a_{i}^{*}})}{1+r}\parallel\nu^{0}-\nu^{*}\parallel<\frac{1}{q}-1$.
Solving for $\parallel\nu^{0}-\nu^{*}\parallel$ and using $1+r>1+r^{\infty}E$ gives the bound used
in Theorem~\ref{Theorem 4}. 

\end{proof}

\noindent\textbf{Acknowledgments:} We thank Jonathan Weare for
introducing JCM to this general area. We also thank the NSF  grant DMS-1613337 and  both SAMSI (DMS-1638521) and Duke TRIPODS (CFF-1934964)
for partial support
of this work.

\bibliographystyle{plain}
\nocite{dinner2017stratification,dinner_mattingly_tempkin_koten_weare_2018,dinner_thiede_koten_weare_2017,cho_meyer_2001,kirkland_neumann_sze_2008,aristoff_bello-rivas_elber_2016,turitsyn_chertkov_vucelja_2011,maier_stein_1996,heymann_2007,heymann_vanden-eijnden_2008,earle_2020,MeynTweedie93,yet11}
\bibliography{paperStratified_gje_jcm_2021}

\end{document}